\documentstyle[11pt]{article}
\def\reals{\hbox{\sl I\kern-.18em R \kern-.3em}}
\def\quats{\hbox{\sl I\kern-.18em H \kern-.3em}}
\def\ints{\hbox{\sl Z\kern-.4em Z \kern-.3em}}
\def\nats{\hbox{\sl I\kern-.16em N \kern.05em}}
\def\rats{\hbox{\sl Q \kern-.83em\vrule height.59em depth0em \kern.87em}}
\def\complexes{\hbox{\sl\kern.50em I\kern-.50em C \kern.05em}}

\def\cC{\cal C}

\def\cF{\cal F}

\def\cR{\cal R}
\def\cS{\cal S}

\input boxedeps.tex
\SetepsfEPSFSpecial
\HideDisplacementBoxes

\begin{document}

\

\centerline{\large{\bf STUDYING SURFACES VIA CLOSED BRAIDS}
\footnote{To appear in J. of Knot Theory and its Ramifications, Volume {\bf 7}, No.
3 (1998)}}

\vspace{.2in}
\centerline{\bf \hskip.2in Joan S. Birman\footnote{partially supported by NSF Grants DMS
91-06584 and 94-02988}  \hskip1in Elizabeth Finkelstein\footnote{partially
supported by a summer grant from Dartmouth College}}\vspace{.1in}

\centerline {{Department of Mathematics   \hskip.6in Department of Mathematics}}
\centerline {{\hskip.1in Columbia University   \hskip1in Hunter College (CUNY)}}
\centerline {{\hskip.1in New York, NY 10027   \hskip1.1in New York, NY 10021 \hskip.1in}}
\centerline {{\hskip.2in jb@math.columbia.edu   \hskip.7in efinkels@shiva.hunter.cuny.edu}}

\

\

\noindent {\bf $\S$ 0. INTRODUCTION.} 

\

In the early $1980's$ Daniel Bennequin
wrote a seminal paper on the theory of contact structures \cite {Bennequin}. The
{\em standard} contact structure $\Delta$ on $\reals^3$ is the kernel of the 
1-form $dz + r^2 d\theta$, where $(r,\theta, z)$ are cylindrical coordinates in
3-space. Thus one attaches to each point in 3-space a planar disc whose position
relative to the coordinate axes is determined by the condition  $dz + r^2
d\theta = 0.$  The question which Bennequin addressed was whether certain other
known contact structures were or were not isomorphic to  $\Delta$. His beautiful
insight was that this question could be referred to a question about braids and
knots.  Bennequin proved that if a knot was everywhere transverse to $\Delta$,
then it could be deformed through transverse knots to a closed braid relative to
the $z$-axis. He then proceeded to study certain incompressible spanning
surfaces bounded by these closed braids, and  using the foliations on those
surfaces which are induced by their intersections with the half-planes $\theta =
\theta_0$, he found invariants of knots transverse to
$\Delta$, under transverse isotopy. Using them he proved that there are
contact structures which are not isomorphic to $\Delta$. 

Several years after Bennequin did his work the first author of this paper  and W.
Menasco initiated, independently, a study of the closed braid  representatives of
a knot or link, using as their chief tool the combinatorics of the very same
foliations on spanning surfaces as had been studied, independently, by
Bennequin. In the course of their work they realized that their ideas were
closely related to the ideas in \cite{Bennequin}, however their interests  were
far from contact structures and moreover they had developed the  applications to
links more fully than Bennequin. They also extended the applications to the
study of surfaces other than spanning surfaces for links. In addition to a
series of papers which Birman and Menasco wrote (full references are given
below) their techniques were exploited in other directions and used in  diverse
ways by J. Los, E. Finkelstein, P. Cromwell and T.S. Fung.  In  particular, the
machinery was applied to the following types of surfaces:

\begin{enumerate}

\item An incompressible spanning surface for a link $L$ (see
\cite{Bennequin}, \cite{Birman-Menasco I},  \cite{Birman-Menasco II},
\cite{Birman-Menasco III}, \cite{Birman-Menasco V}, \cite{Birman-Menasco VI},
\cite{Ko-Lee}, \cite{Los}).  
 
\item A closed incompressible surface in the complement of $L$ (see \cite
{Birman-Menasco IV}, \cite{Birman-Menasco ET},
\cite{Finkelstein}, \cite {Cromwell}).

\item An incompressible 2-sphere in the complement of $L$ which divides 
$L$ into irreducible components \cite{Birman-Menasco IV}

\item A closed 2-sphere which decomposes a non-prime component $K$ of $L$ into
prime summands (see \cite{Birman-Menasco IV}). 

\item A family of immersed `unknotting' discs which $L$ bounds
\cite{Menasco}

\item An immersed ribbon surface which $L$ bounds \cite{Fung}.

\end{enumerate}

While all of the papers cited above used related techniques, those  techniques
were developed in an ad hoc fashion, to fit the needs of the moment. It was 
only after they had been used repeatedly that it was realized that a general
reference might  be in order. That is the purpose of this paper: to review the
Bennequin-Birman-Menasco  machinery, in a systematic fashion,  in the most
general setting which we know at this  time, with an eye toward making it
accessible to the beginning reader. It is our intention that this paper will
serve as a basic reference for the manuscripts
\cite{Finkelstein} and \cite{Fung}. We will also give several new applications
and state some open research problems. The potential  applications do not appear
to have been exhausted.

Here is a guide to the paper. In Section 1 we set up our basic geometry. The
main result in Section 1 is  that incompressible spanning surfaces for
links  and closed incompressible surfaces in link complements admit a singular
foliation which allows  them to be divided into a finite number of foliated
regions, each of which has a canonical embedding in 3-space. See Theorem 1.2. In
Section 2 we describe  ways that we can modify the foliation locally  via isotopy,
and in some cases this involves a change in  the braid representation. In Section
3 we describe certain global properties of the foliation.  In Section 4 we give
examples of ways in which the machinery of Sections 1, 2 and 3 can be applied, and
we pose several interesting open  problems.

The new results in this paper may be summarized as follows:
\begin{itemize}
\item The `change in foliation' move of $\S$2.1 is stated and proved in a new
way, which clarifies its content.
\item The `exchange moves' of $\S$2.4 have been described in a new way which
clarifies the fundamental way in which these moves depend on the surface
foliation.  
\item The graphs which are studied in $\S 3$
are investigated here for the first time as a systematic tool. Theorem 3.1, which
shows how the graphs can be simplified by changes in foliation, exchange moves,
isotopies and deletions of trivial loops, is a new result.
\item The proof of Theorem 4.3, Markov's Theorem without stabilization, in the
special case of the unlink, is new. This concise proof is possible due to the use
of the graphs of $\S 3$.  
\item Example 4.1 is studied here in a
systematic way for the first time. It illustrates that the foliation of a
surface bounded by a knot can indicate how different representatives of the knot
are related. \end{itemize}   

\

\noindent {\bf Acknowledgement:}  We thank the referee for his or her careful
reading of the manuscript and many constructive suggestions. We thank Bill
Menasco for numerous helpful conversations. Many of the results in this paper
are based upon his work,  and his joint work with the first author. We thank
Michael  Hirsch for pointing out an oversight in the statement and proof of
Theorem 2.1. We thank Tahl Nowik for pointing out a missing case in the proof
of Lemma 1.4 and for showing us  that the proof follows from  his work in 
\cite{Nowik}.

\

\

\

\noindent {\bf $\S$ 1. BASIC MACHINERY}. 

\

\

The central idea in the
study of surfaces via closed braids is that when a link is
positioned as a closed braid in $S^3$, there is a foliation of the 
complement  of the
braid axis by discs whose boundary is the axis. This foliation  induces (by
intersection) a  singular foliation of any surface   whose position is 
related to the
position of a given link.  The purpose of this section is to describe the 
general
position arguments which we use to make this foliation amenable to study, 
and 
then to
begin to study it.

{\it We restrict our attention throughout this paper to surfaces which are  embedded
in non-split link complements $S^3 \setminus L$.} We use: 
\begin{itemize}
\item the symbol $\cF$ for an incompressible  spanning
surface for $L$, oriented
so
that the positive normal bundle to each component  has the orientation
induced
by that on $L = \partial \cF$.
\item the symbol $\cC$ for a closed oriented
incompressible surface in $S^3 \setminus L$, and 
\item the symbol $\cS$ when our work refers to both $\cF$ and $\cC$. 
\end{itemize}
We will not consider the case of a 2-sphere which pierces a single 
component of $L$
twice, because that case is rather special, moreover it follows easily 
from the
others. We will not consider split links,  but many of the results hold in
this case as well (see, for example, \cite{Birman-Menasco IV}). We refer the
reader who is interested in  immersed surfaces  to the manuscript of T.S. Fung
\cite{Fung}. That manuscript includes a  discussion of what is known, at this
writing, about the applications of the machinery which  is developed
in this paper to the study of immersions. 

\bigskip

A link type $\cal L$ has a  closed $n$-braid representative
$L$ if there is an unknot $A$ in $S^3 \setminus L$, and a choice of
fibration $H$ of the solid torus $S^3 \setminus A$ by meridian disks, such
that $L$
intersects each fiber of $H$ transversely. Sometimes it will be
convenient to replace $S^3$   by $R^3$
and to think of the fibration $H$ as being by half-planes 
$\{ H_\theta;\theta  \in
[0, 2\pi] \}$  of constant polar angle $\theta$ through the $z$-axis. 
Note that $L$ intersects each
fiber $H_\theta$ in the same number of points, that number being the {\em 
braid
index} $n$ of $L$.  We may always assume that $L$ and $A$ can be oriented so
that $L$ travels around $A$ in the positive direction, using the right hand
rule, as in Figure 1.1.

$$\EPSFbox{1.1.eps}$$
\centerline{\bf Figure 1.1}

\

\noindent It
was proved by Alexander \cite{Alexander} that every link type may be 
represented
 as a closed $n$-braid, for some $n$.  

The braid axis and the fibers of  $H$ will be seen to
serve as a coordinate system in 3-space in which to study closed surfaces  in  $R^3 
\setminus  L$
and spanning surfaces for the link. A singular foliation of each such 
surface
is induced by
its intersection with fibers of $H$. We show (see
Theorem 1.2) that the surface is decomposed into a finite set of canonical
regions, each of which is embedded in 3-space (relative to our
coordinates)
in a canonical way. Later (see Theorem 4.1)  we will show that the surface can be
described via a finite set of combinatorial data associated to this
foliation.
The combinatorial data determines the surface
embeddings by describing how the canonical pieces fit together. 

Our
initial goal in this section is to use some very
simple and well-known general position techniques to place restrictions on
the
foliation of $\cS$. After that we will prove several basic lemmas which
allow
us
to restrict the leaf type and to assume that the foliation has no
$\lq$inessential' leaves. Our first main result is Theorem
1.1. 

\

By general position, we may assume that:

\begin{enumerate}
\item[(i)] The intersections of $A$ and $\cS$ are finite in number and
transverse.
\item[(ii)] There is a  neighborhood $N_A$ of $A$  in $R^3 \setminus L$ such that
each
component of $\cS \cap N_A$ is a disk, and each disk
is radially foliated by its arcs of intersection with fibers of $H$.

\item[(iii)] There is a
neighborhood $N_L$ of $L$ in $R^3$ such that $N_L \cap \cF$ is foliated by
arcs of intersection with fibers of $H$ which are transverse to $L$.

\item[(iv)] All but finitely many fibers $H_\theta$ of $H$ meet $\cS$ 
transversely, and those which do not (the {\it singular fibers}) are each
tangent to $\cS$ at exactly one point in the interior of both $\cS$  and
$H_\theta$.
Moreover,
each point of  tangency is a  local maximum or minimum  or a  saddle
point (with respect to the parameter $\theta$).

\end{enumerate}

\noindent A  {\it singular leaf} in the foliation is one which contains a
point
of tangency with a
fiber of $H$. All other leaves are {\it non-singular}. It follows from
(iv)
that:

\begin{enumerate}
\item[(v)] Each non-singular leaf is an arc or a simple closed curve.
\item[(vi)] Each singular fiber contains exactly one singularity of the
foliation, each of
which is a center or a saddle point.

\end{enumerate}

Note that  each non-singular
leaf in the foliation of $\cC$ 
which is an arc
must have both endpoints on  $A$, since $\cC$ does not intersect $L$.  On
the
other hand, arcs in the
foliation
of $\cF$ may, in principle, have  endpoints on
either $A$ or $L$.
The following lemma, however,  rules out the case of an arc which has both
endpoints on $L$.

\


\noindent {\bf Lemma 1.1.} {\it If $H_\theta$
is a
non-singular fiber of $H$, then an arc in $H_\theta \cap \cF$ cannot have
both of
its endpoints on $L$.}

\

\noindent {\it Proof of Lemma 1.1.} Let $\alpha \in H_\theta \cap \cF$ be an
arc
which has
both of its endpoints on L and let $N_{\alpha}$ be a neighborhood of
$\alpha$
on $\cF$ (see
Figure 1.2). Then $L \cap N_{\alpha}$ has two components, $k$ and
$k^{\prime}$,
which have opposite orientations as subarcs in  the boundary of the
oriented
surface $N_{\alpha}$.
However, $L$ is a closed n-braid, hence it meets each $H_\theta$ in $n$
coherently
oriented
transverse intersections. Since $\alpha$ lies in both $\cF$ and
$H_\theta$, and
since $N_{\alpha}$
intersects $H_\theta$ transversely, this is impossible. $\|$

$$\EPSFbox{1.2.eps}$$
\centerline{\bf Figure 1.2}

\

In view of Lemma 1.1, the non-singular leaves  in the foliation of $\cS$ can
be subdivided into three types:

{\it a}-arcs: arcs which have one endpoint on $A$ and one on $L = \partial
\cS$

{\it b}-arcs: arcs which have both endpoints on $A$. 

{\it c}-circles: simple closed curves

\noindent Each $b$-arc  in a fiber $H_\theta$ separates that fiber into
two
components. Call  the
$b$-arc  {\it essential} if each of these components is pierced at least 
once by
$L$, and {\it inessential} otherwise. Each $c$-circle in a fiber
$H_\theta$
bounds
a disk
${\it D_{c}}$ in $H_\theta$. The $c$-circle  is {\it essential} if $D_{c}$
is  
pierced  at
least once by $L$.  It is {\it inessential} otherwise. We say the
$c$-circle 
is a 
{\it meridian},  if $D_{c}$ intersects $L$ in exactly one point (note that
the
interior of $D_{c}$ may intersect the surface). 

\

\noindent {\bf Lemma 1.2.} {\it All $b$-arcs in the foliation of $\cS$  may
be
assumed to be essential.}

\

\noindent {\it Proof of Lemma 1.2.} If the foliation of $\cS$ contains an
inessential $b$-arc, that arc will cobound with a segment of the axis $A$
a
disk
$\Delta$ in some $H_\theta$ (see Figure 1.3). We can then push the
surface
in along
a neighborhood of $\Delta$ in $R^3$ to remove two points of intersection
of
the
surface with the axis $A$. The removal of the inessential
$b$-arc reduces the number of points at which $A$ intersects $\cS$. Since
the
number of points in $A \cap \cS$  is finite, we can
continue this process until all inessential $b$-arcs have been removed.
$\|$

$$\EPSFbox{1.3.eps}$$
\centerline{\bf Figure 1.3}

\

\noindent {\bf Lemma 1.3.} {\it If the foliation of $\cF$  contains $c$-circles, then
$\cF$ is isotopic to a spanning surface ${\cal F}^{'}$ for $L$ which is foliated
without  $c$-circles.}

\

\noindent {\it Proof of Lemma 1.3.} Suppose that  $c$ is a 
$c$-circle in the foliation of an incompressible spanning surface  $\cF$ for $L$.
Moving through the fibration slightly, we may assume that $c$ is contained in a
non-singular fiber $H_c$. Let
$D_c$ be the disk in $H_c$ bounded by $c$. If $L$ intersects $D_c$,
then there must be arcs in the foliation of $\cF$ of the type ruled out by Lemma 1.1.
Therefore we can assume all $c$-circles in the foliation of $\cF$ are inessential. 

Let $c$ be an inessential $c$-circle in the foliation of $\cF$  contained in
a non-singular fiber $H_c$.  The entire surface cannot be foliated by
inessential circles, since then it would be a compressible torus. Therefore if we
follow  $c$ as it evolves in the fibration, we must arrive in one direction at a
circle   $c_\theta$  containing a singular point $p_\theta$  in a singular fiber
$H_\theta$. The point  $p_\theta$ is
either a center or a saddle point, but if it is a center we can reverse
the flow and find a saddle. So, we may assume it is a saddle. By property (iv)
$p_\theta$ is the only singularity on the fiber $H_\theta$, so the singularity must be as
shown in  Figure 1.4. There are two possibilities, illustrated in Figure 1.4, but
if the second occurs then the first occurs for some other choice of inessential $c$,
so we may assume we are in the situation of case 1.

$$\EPSFbox{1.4.eps}$$
\centerline{\bf Figure 1.4}

\

\noindent Thus $c_\theta$ bounds a disc
$\Delta$ in $H_\theta$ as shown in case 1, and $\Delta$ does not intersect $L$ since
all $c$-circles in the foliation of $\cF$ are  inessential. It may happen that the
interior of the disk $\Delta$ has empty intersection with $\cF$. If so, we may surger $\cF$
along $c_\theta$ as in Figure 1.5 (a). Since $\cF$ is incompressible,  the surface
resulting from this surgery must have one more component than did $\cF$, and  one of the
new components must be a topological $2$-sphere. Discarding the $2$-sphere, we obtain a new
incompressible  spanning surface for $L$. A small isotopy smooths the surface out so that
the saddle $p_\theta$ and the circle $c_\theta$ disappear, and the fibration of the
resulting surface has at least one less saddle singularity  than that of the original
surface $\cF$. 

It remains to consider the case in which the interior of $\Delta$ intersects $\cF$. By
property $(iv)$, there are no singularities of $\cF \cap H_\theta$ in the interior of
$\Delta$, so each component of $\cF \cap int(\Delta)$ must be a 
$c$-circle. Let $c^{'}$ be an innermost such $c$-circle bounding a disk $D$ in
$H_\theta$.  Since $c^{'}$ is inessential, we may surger $\cF$ along $D$, as in Figure
1.5 (b). As before, the surgery produces a closed component $X$ homeomorphic to
$S^2$. This $2$-sphere must be tangent to $H$ at at least two points, so discarding
$X$ eliminates at least two centers. Thus, although the surgery itself  introduces two
new centers, discarding $X$ results in a surface   with no more centers or
saddles than the original.  A finite number of surgeries of this type
result in the disk $\Delta$ having  interior disjoint from $\cF$. A surgery
of the type shown in Figure 1.5 (a) then eliminates it along with its associated
saddle singularity.

Note that the existence of a $c$-circle  implies the existence of a saddle singularity and
that, in each case,  the surgeries  which eliminate it result in a reduction in  the total
number of saddle singularities. Induction on the number of saddle singularities then allows us to
conclude that, after finitely many surgeries,   a surface ${\cal F}^{'}$ is obtained which is
foliated without $c$-circles. Our final observation is that, since the link $L$ is
non-split,  each time we do a surgery the $2$-sphere $X$ must bound a $3$-ball  in
$S^3 \setminus  L$. It follows that each surgery could have been replaced by an isotopy of
the surface through a $3$-ball. Thus, the resulting surface ${\cal F}^{'}$ is isotopic to
$\cF$. $\|$

\

\noindent {\bf Lemma 1.4.} {\it If the foliation of $\cC$  contains inessential
$c$-circles, then $\cC$ is isotopic to a 
surface ${\cal C}^{'}$ in $S^3  \setminus 
L$ which is foliated
without inessential $c$-circles.}

\

\noindent {\it Proof of Lemma 1.4.} The argument which we used to eliminate inessential
$c$-circles in the proof of Lemma 1.3 applies as before, except possibly in the situation
of case 2 of Figure 1.4. In case 2, when the surface is closed and has genus greater than
zero and the $c$-circle $\lq\lq$surgers with itself from the outside," it can happen that
the disk which is bounded by the inner circle after the surgery has non-empty intersection
with the link $L$. We are grateful to Tahl Nowik for pointing this out to us, and for
showing us how this lemma follows from his work in \cite{Nowik}.

In \cite{Nowik}  Nowik studies the intersection of
closed, orientable incompressible surfaces $S$ and  $F$ in an
irreducible orientable $3$-manifold $M$. The   surface $S$ is fixed while  $F$  moves via a
directed isotopy toward a preferred side. Throughout the isotopy $F$ is in 
general position relative to $S$. Now let $M$ be the closure of $S^3 \setminus N$, where
$N = N_A \cup N_L$, $N_A$ is a neighborhood of the axis $A$ and $N_L$ is a neighborhood of
the link $L$.  Let $F$ be the closure of $H_{\theta} \setminus (N \cap H_{\theta})$ and let
$S$ be the closure of $\cC \setminus (N \cap \cC)$. Except for the fact that these surfaces
$S$ and $F$ have non-empty boundary, this is exactly the situation studied in \cite{Nowik}.
Our situation is quite restricted, however, as $\partial F$ consists of meridians of
$\partial N_L$  and a longitude of  $\partial N_A$, and $\partial S$ consists of meridians
of $\partial N_A$. Therefore, since each meridian of $\partial N_A$ has a neigborhood on
$S$ foliated entirely by arcs of intersection with fibers of $H$, all singularities between
$F$ and $S$ which occur as $F$ moves through the fibration are off $\partial F$ and
$\partial S$. The reader can verify that the proof of Theorem 6.1 of \cite{Nowik} holds
in this special case, and that this lemma follows as an immediate consequence. $\|$

$$\EPSFbox{1.5.eps}$$
\centerline{\bf Figure 1.5}

\

\noindent {\bf Lemma 1.5.} {\it All singularities in the foliation of $\cS$
may
be assumed to be saddles.}

\

\noindent {\it Proof of Lemma 1.5.} Suppose there is a center.  Proceeding
through the fibration, we obtain a $c$-circle  in the fiber $H_c$ bounding a
disk $D$
contained in $\cS$ (see Figure 1.6). The circle $c$ is essential, so it bounds a
disk $D_{c} \subset H_{c}$ intersecting $L$ algebraically  non-zero times.
Gluing $D$ to $D_{c}$ along $c$  results in a 2-sphere pierced 
algebraically non-zero times by  $L$. Impossible $\|$

$$\EPSFbox{1.6.eps}$$
\centerline{\bf Figure 1.6}

\

Note that, moving forward through the fibration,  any singular leaf in the
foliation 
is formed by non-singular leaves  moving
together to touch at  a saddle singularity.  After the saddle, the
singular
leaf is
transformed into new non-singular leaves. We call
this process a $\it surgery$ of the non-singular leaves, and label  the
singularities, singular leaves and corresponding surgeries according to
the
non-singular leaves associated to them (see Figure 1.7):

\

{\it type aa:} an $a$-arc surgered with an $a$-arc,

{\it type ab:} an $a$-arc surgered with a $b$-arc,

{\it type bb:} a $b$-arc surgered with a $b$-arc,

{\it type bc:} a $b$-arc surgered with itself or a $c$-circle,

{\it type cc:} a $c$-circle surgered with itself or another $c$-circle.

\

\noindent The preceding lemmas, together with the fact that the foliation
of a closed surface 
contains no $a$-arcs, imply the following theorem.

\

\

\noindent {\bf Theorem 1.1} {\it

\begin{enumerate}

\item[(i)] The surface $\cF$ can be chosen so that every non-singular leaf
in its foliation   is an  $a$-arc or an  essential $b$-arc, and every singularity in its
foliation has  type $aa$, $ab$ or $bb$.

\item[(ii)] The surface $\cC$  can be chosen so that every
non-singular leaf in its foliation is an  essential $b$-arc or an
essential
$c$-circle, and
every singularity in its foliation has  type $bb$, $bc$ or $cc$.

\end{enumerate}}

\

\noindent From now on, we assume that  $\cS$ has been chosen so that its foliation is as
described in  Theorem 1.1.

\

We next show how
the foliation of $\cS$ allows us to decompose $\cS$ into very simple canonical
pieces, which have canonical embeddings in 3-space. Our results will be
summarized, later,  in Theorem 1.2.

Note that no component ${\cal X}$ of $\cS$ is foliated entirely by $b$-arcs, 
since   then ${\cal X}$
is a $2$-sphere, contradicting the incompressibility of $\cS$. Thus, if the foliation of
${\cal X}$ contains no singularities, ${\cal X}$ is either a disk foliated by $a$-arcs (in
the case $\cS = \cF$) or a torus foliated by $c$-circles (in the case $\cS = \cC$).
Otherwise, let $\it U$ be the union of all the singular leaves in the foliation of
${\cal X}$.  Since the fibration
$H$ has finitely many singular fibers, each singular leaf $\lambda$ in $U$
has a

$$\EPSFbox{1.7.eps }$$
\centerline{\bf Figure 1.7}

\newpage

\noindent foliated neighborhood  $N_\lambda$ in $\cS$ such that $N_\lambda \cap U = \lambda$.
If
$\lambda$ has type $aa$, $ab$, or $bb$, $N_\lambda$ 
is one of the foliated
open $2$-cells shown
in Figure 1.8, with the arrows indicating the direction of increasing
$\theta$.

$$\EPSFbox{1.8.eps}$$
\centerline{\bf Figure 1.8}

\

\noindent If
$\lambda$ has type $bc$ or $cc$, $N_\lambda$ is a foliated annulus or
twice-punctured disk, respectively. Figure 1.8 shows this in the case of
the type 1
$bc$ and $cc$ surgeries described in Figure 1.7. The foliated
neighborhoods
for the case
of type 2 surgeries are obtained by reversing the arrows.

The complement of $U$ in ${\cal X}$ is a union $B_1 \cup B_2 \cup...\cup B_k$,
where
each $B_i$ is  foliated entirely by non-singular leaves as shown in Figure
1.9. If we
choose one non-singular leaf in each $B_i$ and declare it
to be a
{\it boundary arc} or {\it boundary circle}, then the union of all the
boundary arcs
and circles determines a {\it decomposition of ${\cal X}$} into {\it regions},
each
of
which is a  foliated neighborhood of one singular leaf.  Let the {\it type}
of
a boundary arc
be {\it type a} or {\it type b} according to whether it is an $a$-arc or
$b$-arc,
respectively.

$$\EPSFbox{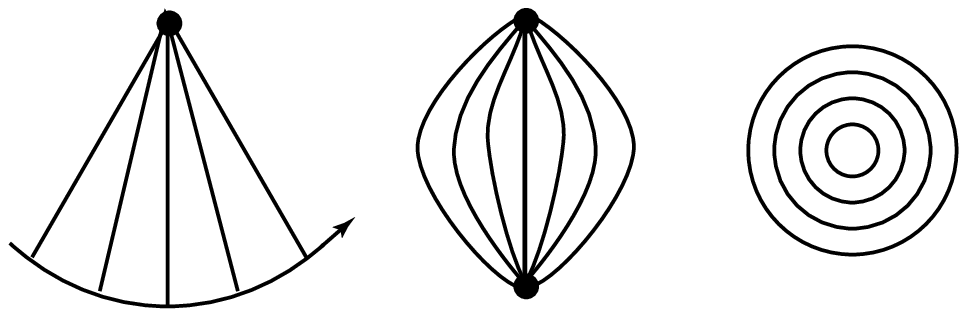 }$$
\centerline{\bf Figure 1.9}

\

From Theorem 1.1, we know that the decomposition of $\cF$ 
contains only  
regions of  type $aa$, $ab$, and $bb$. Since these are $2$-cells, the
decomposition of $\cF$ is a foliated cell-decomposition. We will refer to
these as {\it tiles} of type $aa$, $ab$ and $bb$, and to the
cell-decomposition
a {\it tiling}. 

\

\noindent{\bf Remark 1.1} Observe from Figure 1.8 that an $ab$-tile can be  obtained from a $bb$-tile by 
replacing  one  vertex of the $bb$-tile with a subarc of the link. Similarly, an $aa$-tile can
be obtained from an $ab$-tile by replacing one vertex of the $ab$-tile with a subarc of the link. 

\

Singularities  have a sign. We assume that the fibers $H_\theta$ are
positioned so
that the direction of the positive normal to $H_\theta$ agrees at each
point
with the
direction of increasing $\theta$. At a singular point $s$ in the foliation
the normal to
$\cS$ at $s$ coincides with the normal to $H_\theta$ at $s$. We say that the
singularity at
$s$ is {\it positive} (or {\it negative}) according as the outward-drawn
normal to $\cS$ is
pointing in the direction of increasing (or decreasing) $\theta$. Figure
1.10 shows a
positive singularity viewed in a sequence of fibers of $H$. Let the {\it
type} (respectively {\it sign}) of a region be
equal to the type (respectively sign) of the singularity it contains.

$$\EPSFbox{1.10.eps }$$
\centerline{\bf Figure 1.10}

\

If the foliation of a surface contains $a$-arcs or $b$-arcs, then the axis
$A$ intersects the surface in a finite number of points. Call these points
{\it
vertices} of the tiles.
Each vertex $v$ is an endpoint of finitely many boundary arcs in the
surface
decomposition. Let the {\it type} of $v$ be the cyclic sequence $(x_1,...,
x_r)$,
where  each $x_i$ is either $a$ or $b$, and the  sequence lists the 
types of region boundary arcs meeting at $v$ in the cyclic
order in which they occur in the fibration. The  {\it valence} of $v$ is the
number of
distinct  regions in the surface decomposition intersecting  at $v$. The {\it sign} of $v$
is the
cyclic array
of signs of the regions meeting at $v$.  We define the {\it parity} of $v$ to be  {\it
positive} or {\it negative} according as
the outward-drawn normal to the surface has the same or opposite
orientation as
the braid axis at the vertex. This means that when we view  the positive
side of the surface, the sense of increasing $\theta$ around a vertex will
be
counterclockwise (resp. clockwise) when the vertex is
positive (resp. negative), as illustrated in Figure 1.11.

$$\EPSFbox{1.11.eps }$$
\centerline{\bf Figure 1.11}

\

\noindent {\bf Lemma 1.6} Regions of type $aa,ab,bb$ each have a canonical
embedding in 3-space, which is determined up to an isotopy of 3-space
which preserves the axis $A$ and each fiber of $H$ setwise and up to the
choice
of the sign of the singularity. Regions of type $bc, cc$ each have two
such
canonical
embeddings, which are the same up to a choice of  direction of increasing
$\theta$.

\

\noindent {\it Proof of Lemma 1.6}. We illustrate the proof by showing how
to embed a tile $T$ of type $ab$ in 3-space. Initially, it will be most
convenient to work
in $S^3$, so let us think of the axis $A$ as a circle and the fibers of
$H$
as
discs. Choose an $ab$-tile with vertices at $v_1,
v_2, v_3,$ as in Figure 1.12(a). Let $s$ be the singular point in the
foliation
of $T$ and let $p$ be the $L$-endpoint of the singular leaf in $T$.

$$\EPSFbox{1.12.eps }$$
\centerline{\bf Figure 1.12}

\

The vertices are on $A \cap \cS$, where they have a natural
cyclic order, i.e either $v_1, v_2, v_3$ or its reverse.  One arc of the
singular leaf contained in $T$, say $l_1$, runs from $v_1$ to
$v_3$, separating $H_{\theta_0}$ into 2 halves, one of which contains
$v_2$. The other half of $H_{\theta_0}$ must contain $p$, because the
other
arc of the
singular leaf $l_2$, which runs  from
$v_2$ to $p$, crosses the singular leaf $l_1$. The point of
intersection $l_1 \cap l_2$ is $s$. 

We now pass to $R^3$, choosing the point at infinity so that it
separates $v_1$ and $v_3$. Our fiber $H_{\theta_0}$ is now a
half-plane through the axis, containing the  singular leaf $l_1 \cup l_2$ 
which we just
constructed. We wish to embed the rest of the tile $T$ in 3-space. Since
the
only
singularity on $T$ is the one at $s$, the rest of $T$ will be transverse
to
the fibers it
intersects. The first step is to color the two sides of $T$ light grey (for
positive) and dark grey 
(for negative). Let $w_i$ be a radially foliated neighborhood of $v_i$ on
$T$, i.e. $w_i$
is a sector of a disc. We embed each
$w_i$ in 3-space so that it is orthogonal to $A$ and attached to $A$ at
$v_i$.
Since $v_1$ and $v_3$ are positive and $v_2$ is negative, the light grey sides
of
$w_1$
and $w_3$ face up and the dark grey side of $w_2$ faces up (see Figure
1.12(b)). 

Next choose two overlapping rectangular strips which are neighborhoods
$N_1$
and $N_2$ of the  arcs $l_1$ and $l_2$ on $T$ (see Figure 1.12 (c)) and paste them to
the
$w_i's$ and the singular leaf so that they
are transverse to the fibers of $H$ everywhere except along the singular
leaf. Since the light grey  sides of both $w_1$ and $w_3$ are facing up, we will
need
to give the strip $N_1$ a half-twist as we attach it. (We cannot give it
additional twists because its long edges must be transverse to the fibers
of
$H$.) The half-twist has a sign, determined by the sign of $s$: the light grey 
(or
dark grey) side of $N_1$ faces in the direction of increasing $\theta$
according as
$s$ is positive (or negative). We now pass to $N_2$. The dark grey side of
$w_2$
faces up, so there is a unique way to paste $N_2$ to $w_2$.
The other end of $N_2$ is glued to $L$, and we know how to glue it because
$L$ is
oriented and $S$ has the induced orientation. We also know that $N_2$ is
twisted
just enough so that it can be glued to $N_1$ along the disc $N_1 \cap N_2$
and
so that its long edges are transverse to the fibers of $H$. The full
picture of an embedded $ab$-tile is shown in the top row of  Figure 1.13. 

The type $aa,bb$ cases are similar, and the reader will have no difficulty
in
constructing
the embeddings shown in Figure 1.13, with the help of the pictures given
earlier in
Figure 1.8. The {\it type 1} and {\it type 2} embeddings   for each of
the
$bc$ and $cc$ regions may also be constructed in a similar fashion, using
the
type 1 and
type 2  foliations shown in Figure 1.8, respectively. Note that, in
each case,  the type 1 and type 2  embeddings in $S^3$ are identical, if
the
direction
of increasing $\theta$ and the orientation of $\cS$ are reversed. This
follows
from  
the fact that the corresponding type 1 and type 2  surgeries
are identical, but time-reversed (see Figure 1.7). $\|$

\

\noindent{\bf Remark 1.2} As we mentioned in Remark 1.1, an $ab$-tile can be obtained from a $bb$-tile and
an $aa$-tile can be obtained from an $ab$-tile, by replacing one vertex with a subarc of the link. The
reader can check that the embeddings of the $ab$-tile and $aa$-tile in Figure 1.13 can be obtained in this
manner for a suitable choice of  embedding for the subarc of the link.  We have chosen these
particular embeddings to make the singularities of the foliation visible. We have also distorted the
embedding of the $aa$-tile slightly in order to make its singularity more visible.

\

Summarizing, we have proved:

\

\noindent {\bf Theorem 1.2} 
\begin{enumerate}  
\item [(1.2.1)] {\it Each component of $\cF$ is  a disk foliated by $a$-arcs or  
is  decomposed into canonically embedded
regions of type
$aa,ab,bb$, each of which contains exactly one singular point in the
foliation.} 
\item [(1.2.2)] {\it Each component of $\cC$ is a torus foliated by $c$-circles or is
decomposed into canonically embedded
regions of type $bc,cc$, each of which contains exactly one singular point in the
foliation. }
\end{enumerate}

$$\EPSFbox{1.13.eps }$$
\centerline{\bf Figure 1.13}

\

\ 

\noindent {\bf $\S$2.\  LOCAL CHANGES IN SURFACE DECOMPOSITIONS}.  

\

The surface
decompositions which we have just described are not unique. In this 
section we will describe 
four ways in which they can be changed. In each of the four cases
the possibility of making the change is indicated by examining  
the combinatorics of the tiling. The change is  realized by
an isotopy of the surface which is supported  in a neighborhood 
$N$ of a specified small number of
basic regions of type $aa,ab,bb,bc$ or $cc$, leaving the 
surface decomposition unchanged outside $N$. 

Our goal, of course, is to {\em simplify} the tiling, 
but that depends upon the particular application, and upon having
an appropriate measure of {\em complexity}: 
\begin{enumerate}

\item [(2.1)] The `change in foliation' move of
$\S$2.1  has been used to reduce the valence of particular
vertices in the tiling. See the proofs of Theorem 3.1 and 4.3 
for examples of
how it has been used to reduce complexity. 
  
\item [(2.2)] `Stabilizing along an $ab$-tile', described in $\S$2.2, 
applies to surfaces with boundary.  Let $L=\partial\cF$ be 
a closed $n(L)$-braid. Let $V$ be the number of vertices in 
the tiling of $\cF$. Stabilization along an $ab$-tile will be seen 
to decrease $V$ at the expense of increasing $n(L)$. 
See Theorem 4.2 for an application.

\item [(2.3)] Our  `destabilization' move, described in $\S$2.3, applies to 
surfaces with boundary. It will be familiar to  most readers: 
it removes a `trivial loop' from
the closed braid $L$, reducing both $V$ and $n(L)$. See the Corollary to 
Theorem 4.2 for an application.

\item [(2.4)] Our `exchange move' (which is actually two different but related moves)
is discussed in $\S$ 2.4.  It allows us to detect, by examining the  tiling, how to
`empty pockets' and then remove from the surface. It reduces $V$,
keeping $n(L)$ fixed.
For an application see the proofs of Theorems
3.1 and 4.3.   
\end{enumerate}

\bigskip

\noindent {\bf 2.1 Changes in foliation.} 

\

The choice of a foliation  of $\cS$   is not unique, and our first move
involves ways in which the surface decomposition can be changed by an isotopy of $\cS$ or, equivalently, by an isotopy of
the fibers of $H$ keeping $\cS$ fixed. We may therefore think of the move as either  a change in foliation or a change in
fibration. We have chosen to think of it as a change in foliation. This particular change was introduced in
\cite{Birman-Menasco IV} for $2$-spheres and was modified in \cite{Birman-Menasco V} for certain spanning surfaces.  In
Theorem 2.1, we extend the move  for   arbitrary incompressible spanning  surfaces  and  closed  surfaces $\cS$. As
described in $\S 1$, a singularity  $s$ in the foliation of $\cS$ is  positive (respectively  negative) if the direction of
the positive normal to $\cS$ at $s$ agrees (respectively disagrees) with the direction of the normal to $H_\theta$ which
points in the direction of increasing $\theta$. We label the sign of a singular leaf, surgery, or region, according to the
sign of its associated singularity. Two regions in a surface decomposition are {\it adjacent} if they have a common
boundary arc or boundary circle. 

\

\noindent The reader may wish to look  ahead to Figure 2.7  to see the singular
leaves which correspond to our  changes in foliation, described in Theorem 2.1 below,  for all
possible cases of two regions adjacent at a $b$-arc, when $\cS = \cF$.

\

\noindent {\bf Theorem 2.1.} {\it Suppose that the decomposition of $\cS$ contains regions
$R_1, R_2$ of  the same sign which are adjacent at a $b$-arc or $c$-circle $\gamma$. Then there is
an isotopy $\phi$ taking  $\cS$ to $\cS^{'}$ such that:

\begin{enumerate}
\item [(i)]  $\phi (R_1 \cup R_2) = R^{'}_1 \cup R^{'}_2$,  where $R^{'}_1$ and 
$R^{'}_2$ are adjacent regions of   the same sign in the decomposition of $\cS^{'}$.
\item [(ii)] The decomposition of  $\cS^{'}  \setminus Int(R_1^{'} \cup R_2^{'})$ is the same as 
the decomposition of $\cS \setminus Int(R_1 \cup R_2)$.
\item [(iii)] $\phi(s_1) =  s^{'}_2$ and $\phi(s_2) =  s^{'}_1$, where: $s_i$, $s^{'}_i$ are the
singularities contained in  regions $R_i$, $R^{'}_i$, respectively, for $i = 1,2$ and   $(s_1,
\gamma, s_2)$, $(s^{'}_1, \gamma^{'}, s^{'}_2)$ are their cyclic orders in the fibration relative
to  the common boundary $\gamma$, $\gamma^{'}$ of their respective regions.    
\item [(iv)] If $\cS = \cF$ (so by Theorem 1.2.1 $\gamma$ is
a $b$-arc), $\phi$ decreases the valence of each vertex $v$ and $w$ in $\partial\gamma$ by one.
\end{enumerate}}

\

\noindent {\it Proof of Theorem 2.1.} In Figures 2.1 (a) and (c) we have
illustrated two different embeddings  of a subdisk of  $R_1 \cup R_2$ in $3$-space.
The fibers of $H$ are viewed as horizontal planes, and the polar angle function is the
height function. Note that each pair of singularities  in  (a), (c) 
bounds an arc $\mu$, $\mu^{'}$, respectively,  in $\cS$  which is transverse to
the fibration.  To visualize the isotopy of $\cS$ which takes the embedded disk shown
in (a) to the  disk shown in (c), pull $\cS$ along as  the arc $\mu$ is isotoped to the
arc $\mu^{'}$.  Note that the singularities $s_1$, $s_2$ which are endpoints of $\mu$
have heights opposite to those of their images $s^{'}_2$, $s_1^{'}$, respectively,  
under the isotopy.  Figure 2.1 (b) shows the intermediate stage of  the isotopy at
which the endpoints of the image of $\mu$ have the same height. The foliations which
correspond to the embeddings  of Figure 2.1  are shown directly below them in
Figure 2.2. We define a {\it change in foliation} as the passage from the embedding in
(a) either to the embedding in (c) or to a  $120^{0}$ rotation of the embedding in (c) about
the direction of increasing $\theta$, according to the left hand rule. We will  prove that
under the given hypotheses there is no obstruction to such a change and that, in the case of
regions in the decomposition of a spanning surface  adjacent at a $b$-arc $\gamma$, the change 
decreases  the valence of each vertex of $\gamma$ by one. 

Note that the embedding in Figure 2.1(c) is obtained from the embedding in Figure 2.1(a) by a
$120^{0}$ rotation about the direction of increasing $\theta$, according to the left hand rule. 
Therefore, to show that  there is no obstruction to a change
in foliation it suffices, by symmetry, to show that there is no obstruction to the passage from
the embedding in Figure 2.1(a) to the one in Figure 2.1(c).

$$\EPSFbox{2.1.eps }$$
\centerline{\bf Figure 2.1}

\

\noindent It will be convenient to
work in $S^3$, so assume that each fiber $H_\theta$ is a meridian disk. Up to a choice
of orientation of $\cS$,  the surgeries associated to the singularities $s_1$, $s_2$ 
are shown in Figure 2.3, as we proceed forward through the fibration.  Thus, there
is a disk $\Delta \subset R_1 \cup R_2$ with  embedding  as shown in
Figure 2.1 (a) and   foliation as shown in Figure 2.2 (a).  
Note that the arc  $\mu \subset \Delta$ is contained in  the disk $D$ shown in Figure
2.4 (a).  If $Int(D)$$\cap S=${\em Int}$(\mu)$  and $D \cap L = \emptyset$,  then we can
isotope $\Delta$ to the  disk $\Delta^{'}$ of Figure 2.1 (c) by pulling $\cS$ along
in a small neighborhood of $D$ as we push $\mu$ through $D$ to $\mu^{'}$ as shown in
Figure 2.4. Since  this isotopy does not change $\cS  \setminus$$\Delta$, it leaves
the decomposition of  $\cS   \setminus$$(R_1 \cup R_2)$ unchanged as well.  Also,
since  it is the identity on $\cS  \setminus$$\Delta$ and  since  $\Delta^{'}$ has the 
foliation shown in Figure 2.2 (c), it is a simple matter to check  that it takes
$R_1 \cup R_2$ to $R^{'}_1 \cup R^{'}_2$,  where $R^{'}_1$ and  $R^{'}_2$ are adjacent regions
of the same sign in the decomposition of $\cS^{'}$. Thus, to prove (i),(ii),(iii) of the
theorem, we only need  to show that we can first isotope $\cS$, without changing its
decomposition, so that such a disk $D$ with interior disjoint from $\cS \setminus \mu$ and $L$
exists. This isotopy followed by the isotopy from  Figure 2.1  (a) to Figure 2.1 (c) is exactly
the isotopy $\phi$ described in (i)(ii)(iii) of the theorem.

Note that $\mu$ separates  the desired disk  $D$ into two subdisks which are
identical up  to a choice of increasing $\theta$. By symmetry, 
we need only show that one of the subdisks has interior disjoint from
$\cS$ and $L$. It therefore suffices to find  an arc $\alpha$ as shown
in Figures 2.3 and 2.4 which  has interior
disjoint from $\cS$ and $L$, as we proceed through the fibration between the fibers
containing $s_1$ and $s_2$.

$$\EPSFbox{2.2.eps }$$
\centerline{\bf Figure 2.2}

\

$$\EPSFbox{2.3.eps }$$
\centerline{\bf Figure 2.3}

\

$$\EPSFbox{2.4.eps }$$
\centerline{\bf Figure 2.4}

\

We first prove the case in which    $R_{1}$ and $R_{2}$  are adjacent at a
$c$-circle $c$.  Up to a choice of orientation of $\cS$, we may assume the singularities $s_1$,
$s_2$ in $R_1$, $R_2$ are positive.  Let $H_{c}$ be the fiber
containing $c$,   let $D_{c}$ be the disk in $H_{c}$ bounded  by $c$, and  let $A_{c}$ be the
annulus $H_{c}  \setminus D_{c}$.   Let 
$H_1$ and  $H_2$  be the fibers occurring just after $s_{1}$ and $s_{2}$ in the
fibration, respectively. We shall not distinguish between  $c$ and any circle isotopic to $c$ on
$\cS$ between   $H_1$ and $H_2$.  Since $R_1$ and $R_2$ are adjacent at a $c$-circle,
$R_1$ has type $bc$ or type $cc$. 

Suppose $R_1$  has type $cc$.  It has
either a type 1 or type 2 embedding in $S^3  \setminus L$   (see Figure 1.13), and
Figure 2.5 (a) and (b) illustrates the corresponding surgeries resulting in the
positive singularity $s_1$.  Suppose that  $R_1$ has a  type 1 embedding. Since $s_1$
occurs just before $H_1$,  there is an arc $\alpha \subset D_{c}$ in $H_1$ as shown
in  Figure 2.5 (a)  with interior disjoint from    $L$ and $\cS$. However the link
must (and the surface may) intersect $D_{c}  \setminus \alpha$ at this point in the 
fibration.  We must show that we can isotope $\cS$, without changing its
decomposition, so that the interior of $\alpha$ does not intersect $\cS$ or $L$, as we
proceed forward through the fibration between $H_1$ and $H_2$.

$$\EPSFbox{2.5.eps }$$
\centerline{\bf Figure 2.5}

\

Let $X$ be the 3-ball which is formed by the union of  all the disks $D_{c}$ contained
in the fibers  between $H_1$ and $H_2$. If   no braiding of $L$ or surgeries occur in
$X$,   then there is nothing to  obstruct $\alpha$ as we proceed forward through the
fibration  between $H_1$ and $H_2$. Suppose that braiding of $L$ or  singularities  
occur in $X$.  Since none of the  singularities involve $c$ and since $L \cap c =
\emptyset$, we can isotope the intersection of the interior of $X$ with $\cS \cup
L$    forward  in a direction transverse to the fibration, until all  braiding and
singularities have been pushed out of $X$. Since this isotopy moves  $\cS$ and $L$ in
a direction  transverse to each fiber,  $L$ remains transverse to the  fibration, and
the decomposition of that part of $\cS$ affected by the isotopy is  not disturbed. 
Afterwards,  $\alpha$ can remain fixed  as we proceed forward in the fibration through
$X$, with its interior disjoint from  $L$ and $\cS$. Then, since  any surgery
resulting in a positive $s_2$ must  be represented by the dotted line in $A_c$ shown
in  Figure 2.5 (a), the singularities and braiding of $L$ which we just pushed
forward out of $X$ can be assumed to occur after $s_2$ in the fibration. In
addition,   the singularity  $s_2$ does not cause $\cS$ to intersect the interior of 
$\alpha$. Thus,  we obtain the desired subdisk of $D$ (foliated by $\alpha$) with
interior disjoint from $\cS$ and $L$. Note that this process is only possible when
$s_1$ and  $s_2$ have the same sign. If $s_2$ has sign  opposite to that  of $s_1$, an
attempt to push the braiding of $L$ and singularities of $\cS$ forward past $s_2$ can
result in a link which  is not transverse to the fibration and a drastic change in the
decomposition of $\cS$.

Suppose that $R_1$ has a type 2 embedding. There is an arc
$\alpha \subset A_{c}$ in $H_1$ as shown in  Figure 2.5 (b)  with interior disjoint
from    $L$ and $\cS$, and   any surgery
representing a positive $s_2$ is represented by a dotted line in $D_c$. Let $X$ be the thickened
cylinder formed by the union of all annuli $A_c$ between $H_1$ and $H_2$. As in the
type 1 case, none of the  singularities in $X$ involve $c$ and  $L \cap c = \emptyset$.
Therefore we can isotope the intersection of the interior of $X$ with $\cS \cup L$   
forward  in a direction transverse to the fibration, until all  braiding and
singularities have been pushed out of $X$, and there is nothing to obstruct $\alpha$
between $H_1$ and $H_2$. Then, since  the surgery corresponding to  $s_2$ must be
represented by a dotted line in $D_c$, we can push the braiding of $L$ and
singularities of $\cS$ forward in the fibration past $s_2$. As in the preceding case,
the singularity  $s_2$ does not cause $\cS$ to intersect $\alpha$, and we again 
obtain the desired subdisk of $D$ with interior disjoint from $\cS$ and $L$. 

If $R_1$ has type $bc$, then it must have a type 1 embedding, since otherwise it
is not adjacent to $R_2$ at a $c$-circle. The desired subdisk of $D$ is
obtained in this case by essentially the same proof as for the case in which $R_1$ has
type $cc$ and a type 2 embedding. This
completes the proof of the theorem in the case in which  $R_1$ and $R_2$ are adjacent
at  a $c$-circle.

The proof of (i)(ii)(iii) in the case in which  $R_1$ and $R_2$  are adjacent at a  $b$-arc is
similar in flavor to the previous case.  Suppose that  $R_1$ and $R_2$  are adjacent at a  $b$-arc
$\beta$ with vertices $v$ and $w$. In this case, we have the situation shown in Figure 2.6, as we
proceed forward in the fibration. As in the previous case, let $H_1$ and  $H_2$ be the fibers which
occur just after the singularities $s_1$ and  $s_2$, respectively. Again, there is an arc $\alpha$ in
$H_1$ which does not intersect $\cS$ or $L$ (see Figure  2.6). We must show that
$\alpha$ does not intersect $\cS$ or $L$, as we proceed forward through the fibration from $H_1$ to
$H_2$. Note that the arc $\beta$ separates $H_1$ into two disks. Let $D_{\alpha}$ be the disk 
which contains $\alpha$ (it  is shaded in Figure  2.6).  Since $\beta$ does not
surger between $H_1$ and  $H_2$,  there is such a disk in each fiber between $H_1$ and
$H_2$. Let $X$ be the $3$-ball which is formed by the union  of all the  disks
$D_{\alpha}$. If no singularities or braiding of $L$ occur in $X$, then $\alpha$ can
remain fixed in $D_{\alpha}$  with interior disjoint from  $\cS$ and $L$ as we proceed
forward in the fibration through $X$. Since $\beta \cap L = \emptyset$ and since none
of the  singularities which may occur between $H_1$ and  $H_2$ involve $\beta$, we can
push all singularities and braiding of $L$ forward in the fibration out of the
$3$-ball $X$ as in the case when the regions are adjacent at a $c$-circle. Now note
that the surgery corresponding to $s_2$ is represented by a dotted arc in $H_1
\setminus D_{\alpha}$. Thus, we can push the braiding of $L$ and singularities of
$\cS$ forward in the fibration past $H_2$, and  $s_2$ does not cause $\cS$ to
intersect $\alpha$ or force $L$ to a position that is not transverse to the fibration.
As in the $c$-circle case, the desired subdisk of  $D$ with interior disjoint from
$\cS$ and $L$ is foliated by the arcs $\alpha$.  This completes the proof of
(i)(ii)(iii) of the theorem.

$$\EPSFbox{2.6.eps }$$
\centerline{\bf Figure 2.6}

\

Our proof of (iv)  requires a careful inspection of the region foliations before and after the
isotopy $\phi$.  Note that, since  $\cS = \cF$ is a spanning surface, any two  regions adjacent
at a $b$-arc $\gamma$ must have type $bb$ or $ab$. This yields three cases, up to symmetry, for
the types of $R_1$ and $R_2$: (1) both  are $bb$-tiles; (2) $R_1$ is a $bb$-tile
and $R_2$ is an $ab$-tile or (3) both are $ab$-tiles. We only need to show that the
singular leaves corresponding to the change in foliations for these three cases are  as shown in
(1), (2) and (3) of Figure 2.7, respectively, since each of these figures  illustrates   the
required reduction in valence for the vertices $v$ and $w$ of $\gamma$. The first column after
the arrows in Figure 2.7 illustrates the result of a change in foliation as described in Figure
2.1. The second column after the arrows shows  the result of the second type of change in
foliation, i.e.,   the passage from Figure 2.1(a) to a $120^{0}$ rotation  of Figure
2.1(c). Since the proofs are identical, we will only show that the  singular leaves in the first
column after the arrows are the result of the change in foliation shown in Figure 2.1.

$$\EPSFbox{2.7.eps }$$
\centerline{\bf Figure 2.7}

\

Suppose we are in case (1). Thus  $R_1$
and $R_2$  are two  $bb$-tiles adjacent  at a $b$-arc  with vertices $v$ and $w$. Note
that  the  arcs labeled $1, 2, 3, 4, 5,6$ in the boundary of $\Delta$ in Figure 2.2 (a) are
portions of $b$-arcs, and  the unlabeled arcs  are transverse to the foliation. Therefore, in
order for the leaves in the foliation to match up,  $\Delta$ must be situated  in $R_1 \cup R_2$
as shown in Figure 2.8 (a).  Since $\phi$ does not change  the foliation of $R_1 \cup R_2$ 
outside $\Delta$, the foliation of $R_1^{'} \cup R_2^{'}$ can be obtained by replacing the
foliation of $\Delta$ with that of $\Delta^{'}$. Thus, $R_1^{'} \cup R_2^{'}$ is foliated as
shown in Figure 2.8 (b). The singular leaves corresponding to a change in foliation in case
(1) are therefore as shown in Figure 2.7 (1). This concludes the proof of (iv) in case
(1). The proof for  cases (2) and (3) is identical. Using Remark 1.1, the reader will be able
to see immediately that the foliation changes resulting from $\phi$ in these cases are as shown
in Figures  2.7 (2) and (3), respectively. This concludes the proof of the theorem. $\|$

$$\EPSFbox{2.8.eps }$$
\centerline{\bf Figure 2.8}

\ 

\noindent{\bf Remark 2.1} The reader can  use the  same method which is used in the proof of
Theorem 2.1 (iv) to check  that if a vertex $w$ has type $(a,b)$ and sign
$(\pm,\pm)$ then, after a change in foliation, $w$ has type $(a)$. The singular leaves
corresponding to this change in foliation are shown in Figure 2.9.

$$\EPSFbox{2.9.eps }$$
\centerline{\bf Figure 2.9}

\

\

\noindent {\bf 2.2 Stabilization along $ab$-tiles.}  

\

In this subsection and the next we
will describe two modifications for an incompressible spanning surface $\cF$   which change
the braid representative $L$ and the tiling, but do not change the link type $\cal L$. 
Both of these moves change the  braid index. They were, in effect, introduced by 
Bennequin in \cite{Bennequin}. The $\lq\lq$stabilization move" which we describe in this
subsection adds a trivial loop to the braid representative $L$; and the $\lq\lq$destabilization 
move" which we consider in the next subsection deletes one. In this respect, the two
moves are mutually inverse with respect to the  closed braid structure on $\cal L$.
However, they are not mutually inverse with  respect to their effect on the tiling of
$\cF$. Indeed, the stabilization move results in the elimination of a negative vertex 
whereas the inverse of the destabilization  move results in the addition of a positive
vertex.  The stabilization move deletes an $ab$-tile, and  changes nearby tile types. The
inverse of the destabilization move, adds an $aa$-tile, and leaves the rest of the
tiling  unchanged.

The modification in  the  tiling of $\cF$ which we
call  {\it stabilization along an $ab$-tile}, is shown in Figure 2.10(a).  It  is 
realized by  pushing a subarc $\alpha$ of $L$ along a disc neighborhood of part of the singular
leaf in an $ab$-tile $T$. The boundary of the neighborhood may be chosen to be
transverse to leaves of  the foliation and to the fibration, so the
result is a new closed braid representative $L^{'}$ of the link,
bounding a new surface $\cF^{'}$ which is  tiled, but with one less tile than
$\cF$ and one less negative  vertex. In fact, letting $N$ be  the union of all tiles
of $\cF$ with $v$ as a vertex, we can summarize the effect of
stabilizing along the $ab$-tile $T$ as follows:

\begin{enumerate}

\item[(1)] $T$ and  its negative vertex are deleted.

\item[(2)] Any $ab$-tile in $N \setminus T$  is
transformed into an $aa$-tile.

\end{enumerate}

$$\EPSFbox{2.10.eps }$$
\centerline{\bf Figure 2.10}

\begin{enumerate}

\item[(3)] Any $bb$-tile in $N$ is transformed into an $ab$-tile.

\item[(4)] The tiling of $\cF^{'}$ is the same as the tiling of $\cF$ outside $N$.

\item[(5)] The resulting braid representative $L^{'}$ is obtained from $L$ by the
addition of a trivial loop, as shown in Figure 2.10(c). This increases the braid index
by one.

\end{enumerate}

\noindent To see (1)(2)(3)(4), examine Figure 2.10(b), which   shows one  sample $N$,
flattened out. The change in braid representative described in (5) can be seen by
noting that Figure 2.10(c) is just a projection   parallel to the
axis of the link   described in  Figure  2.10(a).

\

\

\noindent{\bf 2.3 Destabilization via a type $(a)$ vertex.} 

\

We will now  show that
whenever there is a vertex of valence one in the tiling of $\cF$, then  it has type
$(a)$, and we can modify the tiling in the manner shown in  Figures 2.11 (a) and (b).
We call this modification   {\it destabilization via  a type $(a)$ vertex}.

$$\EPSFbox{2.11.eps }$$
\centerline{\bf Figure 2.11}

\

A vertex $v$ of  valence 1 in the tiling occurs when two of the edges in a single tile
$T$  are identified in $\cF$. Such a vertex must have type $(a)$. For if it  has type
$(b)$, the only other possibility, then there is a type $bb$ (respectively $ab$) tile in
the tiling of $\cF$ which has two adjacent $b$-edges identified. This requires that $T$
have only 3 (respectively 2) distinct vertices instead of 4 (respectively 3). The tile
$T$ contains exactly one singular point, and  this singularity occurs when two
non-singular leaves come together. This means that just before the singularity there
were two components of $H_{\theta} \cap \cF$ emerging from a single point of $A \cap \cF$,
which is impossible. Thus $v$ has type $(a)$ and, since the only tile   which has
a pair of adjacent $a$-edges meeting at a common vertex is an $aa$-tile, it follows
that  $T$ is  $aa$-tile. 

We know there is exactly one canonical embedding for
an $aa$-tile in 3-space, shown in Figure 1.13.
Identifying two edges of an  $aa$-tile with endpoints on a common vertex therefore 
yields  the canonical embedding in 3-space for $T$   shown in Figure 2.11 (a). Notice
that there is a radially foliated disk $D$ in $T$ cut off from $\cF$ by
the arc of the singular leaf of $T$ with both endpoints on $L$.
Therefore, there must be a trivial loop in the braid representation of
$L$.  Destabilization via the type $(a)$ vertex $v$ is realized by pushing $L$ through $D$ 
to the singular leaf in $T$, and then modifying this new link representative so that it
becomes transverse to all fibers of $H$ (see Figure 2.11 (b)). The result is a new closed
braid representative
$L^{'}$   bounding a new surface $\cF^{'}$, with changes summarized as follows:

\begin{enumerate}

\item[(1)] The $aa$-tile $T$ and its type $(a)$ vertex $v$ are deleted.

\item[(2)] The tiling of $\cF^{'}$ is the same as the tiling of $\cF$  outside $T$.

\item[(3)]  $L^{'}$ is obtained from $L$ by deleting a trivial loop, which decreases  the
braid index by one.

\end{enumerate}

\

\noindent {\bf 2.4 Exchange Moves.}  

\

The modification which was 
described in
$\S$2.3 was based upon the existence of a vertex of valence 1 
in the tiling.
In this subsection we consider modifications which are based upon 
the existence of a
vertex of valence 2.  

The precise definition of our exchange move (given in Theorem 2.2 below) 
is based upon the tiling of our surface, however the reader may find it
easier to understand the move by looking first at its effect on a knot 
diagram.   See Figure 2.12 for an example. The labels $n_j$ on the
strands are weights, indicating $n_j$ parallel strands. We allow any 
type of braiding
on the $n_1 + n_2$ (respectively  $n_1 + n_3$) strands in the boxes 
which are labeled
$X$ (respectively $Y$). Assume that the braid axis $A$ is the $z$-axis, 
and that the
arc which is labeled $n_3$ lies in the $xy$-plane. Up to isotopy of 
$S^3$, an
exchange involves an isotopy of $L$ which moves the arc which 
is labeled
$n_2$ from a position which is a little bit above (or below) the 
$xy$-plane to a
position which is a little bit below (or above) the $xy$-plane.  
However, as we
shall see, this is only part of the story, because all of our moves are guided by
the foliation of the surface, and exchange moves are too.

We shall define exchange moves in two different situations 
(parts (2.2.1) and (2.2.2) of Theorem 2.2). The phenomena  will seem, initially,  to
be very different in the two cases, but after the proofs are complete we will be able
to understand why they are in fact very closely related.

$$\EPSFbox{2.12.eps }$$
\centerline{\bf Figure 2.12}

\

\noindent {\bf Theorem 2.2 \cite{Birman-Menasco IV}, \cite{Birman-Menasco V}}
\begin{enumerate} 

\item [(2.2.1)]  {\it Suppose that
$\cF$ has a vertex $v$ of valence 2, type $(a,b)$ and sign 
$(\pm,\mp)$.  Then
$L=\partial\cF$ admits an isotopy to a new closed braid 
$L^\prime$, such that:}
\begin{enumerate}
\item [(a)] {\it The isotopy from $L$ to $L^\prime$ is a push 
across a disc
$\Delta\subset\cF$, where $\Delta$ is a neighborhood of a disc 
$\delta$ which
is contained in the union of two $ab$  tiles in the decomposition 
of $\cF$.}

\item [(b)] {\it The decompositions of $\cF$ and $\cF^\prime = 
\cF\setminus
\Delta$ are identical everywhere except on $\Delta$.  
The tiling of $\cF^\prime$ has two fewer vertices than the 
tiling of $\cF$.
The braid indices of $L$ and $L^\prime$ coincide.  }

\item [(c)] {\it The isotopy $L \to L^\prime$ modifies $L$ near $A$ as
depicted in Figure 2.12.} 
\end{enumerate}

\item [(2.2.2)] {\it Suppose that $\cS$ 
has a vertex $v$ of
valence 2, type $(b,b)$ and sign $(\pm,\mp)$.  Then $L$ admits an 
isotopy to a
new closed braid $L^\prime$ and $\cS$ admits an isotopy to a
new surface $\cS^\prime$ such that:}
\begin{enumerate}
\item [(a)] {\it The decompositions of $\cS$
and $\cS^\prime$ are identical except on the two given $bb$ tiles in
$\cS$. The tiling of
$\cS^\prime$ has two fewer vertices and two fewer singular points 
than the tiling of
$\cS$.   The braid indices of $L$ and $L^\prime$ coincide.} 

\item [(b)] {\it The isotopy $L \to L^\prime$ modifies $L$  near $A$ as
 as depicted in Figure 2.12.} 
\end{enumerate}
\end{enumerate}

\noindent {\it Proof of Theorem 2.2, part (2.2.1):} (See page 597 of
\cite{Birman-Menasco V})  A vertex of  type $(a,b)$ occurs
when two $ab$-tiles $T_1$ and $T_2$ are identified along 
corresponding $a$ and
$b$ edges, as in Figure 2.13(a). The vertex $v$ is joined to 
$\partial \cF$ by the 
common $a$-edge in $T_1$ and $T_2$ and to a second common vertex $w$   by the common
$b$-edge.

Let $\delta$ be
the disc in $T_1\cup T_2$ which is cut off by the portions of the singular 
leaves which connect  $w$ to $L$ in
$T_1\cup T_2$. We know  the foliations of our two tiles from 
Figure 1.8, and we
know that a neighborhood of $w$  is radially foliated. We 
may then find an arc
$\alpha^{\prime}$ (see Figure  2.13(a)) which begins and ends 
on $L$ and is everywhere
transverse to leaves of the foliation, which  cobounds with 
a subarc $\alpha$ of $L$ 
a neighborhood $\Delta$ of $\delta$ on $\cF$.  The disc $\Delta$ is shaded in
lightly in Figure  2.13(a). Push $\alpha$ through
$\Delta$ to $\alpha^\prime$, crossing the axis twice in the process, 
to obtain a new
closed braid representative  $L^{'}$ of $\cal L$ in $S^3$. The new 
 representative
$L^{'}$  bounds a new surface $\cF^{\prime}$, isotopic to $\cF$. 
The new surface
$\cF^{\prime}$ is also tiled, but has two fewer tiles and two 
fewer vertices than the
tiling of $\cF$. In fact, if $\cR$ is the union of our two 
$ab$-tiles together with
all tiles in the decomposition of $\cF$ having  $w$ as a vertex, 
then the
decomposition of $\cF^{'}$ is the same as the decomposition of 
$\cF$ outside $\cR$. An
examination of Figure  2.13(a) shows that the changes in the decomposition 
on $\cR$ are as
follows:  our two $ab$ tiles are deleted; (2) any other $ab$-tile 
in $\cR$ becomes
an $aa$-tile and (3) any $bb$-tile in $\cR$ becomes an $ab$-tile.

$$\EPSFbox{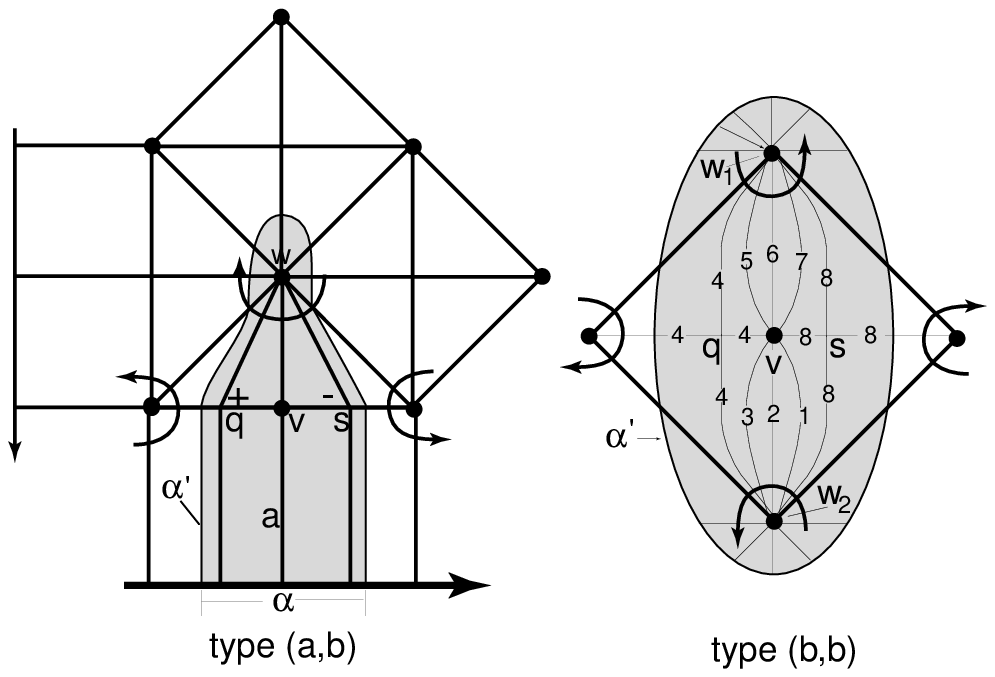 }$$
\centerline{\bf Figure 2.13}

\

Note that
the sense of rotation of the flow is in opposite directions about the 
vertices $v$ and
$w$. Thus, the algebraic rotation number  of the 
flow about the
deleted vertices is zero, so the algebraic rotation number of 
$L^{'}$ about $A$ is the
same as that of $L$, namely $n$.  From this it follows that 
the braid index of $L^{'}$ 
is $n$, because $L^{'}$ is transverse to every fiber of $H$. 

It only remains to show that $L$ and $L^{'}$ are related by an 
exchange move. 
With this in mind, we study the embedding of $\Delta$ in $\reals^3$.
Since $v$ has sign $(\pm,\mp)$,  the two  singularities in $\Delta$ 
have opposite
signs. Thus, since the flow is   in opposite directions about the 
vertices $v$ and
$w$, it follows that $\Delta$ must be embedded as shown in Figure 2.14(a). 
 Now since all
$b$-arcs in the foliation are essential, the link itself must be an 
obstruction to
their removal, as we have shown  in Figure 2.14(a). We can push all 
braiding away from
$A$, into the box labeled $X$. Projection parallel to $A$ now clearly 
shows that when
we push $\alpha$ through $\Delta$ to $\alpha^{\prime}$, we are doing 
an exchange 
move. In this case $n_3 = 1$.

$$\EPSFbox{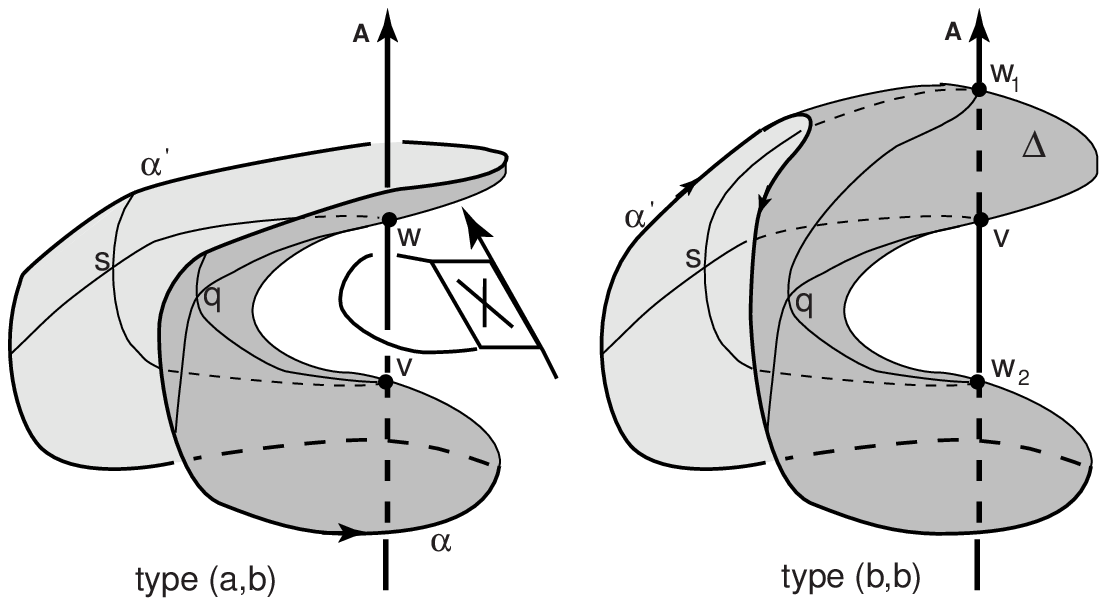 }$$
\centerline{\bf Figure 2.14}

\

\noindent {\it Proof of (2.2.2):} (See page 130 of
\cite{Birman-Menasco IV}). We now suppose that $v$ has sign $(\pm,\mp)$ 
and  type $(b,
b)$.  The geometry in (2.2.2) is in many ways similar to that in 
(2.2.1), because an $ab$ tile may be thought of as having been 
obtained from a $bb$
tile by deleting a neighborhood of a vertex. However,  in (2.2.2) the region 
near $v$ is in the
interior of ${\cal F}$, whereas we will be modifying $L$. The modifications are
suggested by Figure 2.15.

In order to justify this picture and prove (2.2.2) we will
need to:  
\begin{enumerate}
\item Show that $\cS$ and $L$ are positioned as in the left picture in Figure 2.15; 
\item Find structures in the foliation of the
two $bb$ tiles which meet at $v$ which will tell us how to modify $L$ to the
middle picture in Figure 2.15. 
We already know
that $ab$ tiles are like $bb$ tiles with a deleted vertex, and the structures will
appear out of a ``shadow" of our two $bb$ tiles;  
\item Show that the decomposition of $\cS$ is left invariant under the modification,
even when $\cS$ is a spanning surface and $L$ its boundary;  
\item Show that, on the other hand, the order in which tiles on the surface are 
intersected by $A$ {\it is} changed;  
\item Show that the change in the braid
structure of $L$ is an exchange move, as illustrated in Figure 2.12, and finally; 
\item Show
that after the change we have inessential $b$ arcs which can be removed by 
an isotopy of ${\cal S}$ to $\cS^\prime$ in the right picture in Figure 2.15.
\end{enumerate}

$$\EPSFbox{2.15.eps }$$
\centerline{\bf Figure 2.15}

\  

Let $T_1$ and $T_2$ be our two $bb$ tiles and let $q$ and $s$ be the singular
points of the two regions which meet at $v$. Let $w_1$ and $w_2$ be the 
other vertices of
the two edges which meet at $v$. See Figure  2.13(b), which will seem easy
because of the analogy with Figure  2.13(a). 
We have only made one small change in the passage
from Figure  2.13(a) to  2.13(b), i.e. we have changed the signs of the vertices to
stress the fact that the given vertex of type $(b,b)$ and sign $(\pm,\mp)$ could
be either positive or negative, whereas in the case of a vertex of type $(a,b)$,
i.e. as in Figure  2.13(a), it is always positive.  
Portions of the singular leaves which connect $w_1$ to $w_2$ 
fit together
to cut off a disc $\delta$ on $\cS$, and we have
sketched in a neighborhood $\Delta$ of $\delta$ on $\cS$, chosen so that $\partial
\Delta$ is everywhere transverse to the foliation.  
Notice that the gradient flow on $\cS$ is always oriented in opposite senses
around the two endpoints of each leaf which joins two points of $A \cap \cS$. This
implies that $\partial \Delta$ has algebraic rotation number 1 about the axis $A$.
Since $\partial \Delta$ is everywhere transverse to the foliation, it then follows
that $\partial \Delta$ is a 1-braid. The axis $A$ pierces $\Delta$
from alternating sides at $v$, $w_1$ and $w_2$.
Using this and the fact that $s$ and $q$ have opposite signs, the reader should have
no trouble using the procedure described in $\S$1 to verify that  $\Delta$ must be
embedded as shown in  Figure 2.14(b).  

We will need to study this embedding in some detail, and so
it will be convenient to assign numbers $0 \leq \theta_1 < \theta_2 < \dots <
\theta_8 < 2\pi$ to the leaves in the foliation of $\Delta$, as in Figure  2.13(b),
where the label $``i"$ means $\theta_i$.  We need to use these leaves to create a
region in 3-space in which to describe  the exchange move and the isotopy of $\cS$ in
a controlled fashion.

Each non-singular $H_{\theta}$ meets $\cS$ in a unique arc $b( \theta )$ which 
joins $v$
to $w_1$ or $w_2$ and cobounds with a part of the axis $A$  two discs in the
fiber $H_{\theta}$. If $b( \theta )$ joins $v$ to $w_1$ (resp. $w_2$), let
$\mu_{\theta}$ be the disc which does not contain $w_2$ (resp. $w_1$). We have 
sketched in
one such arc $b( \theta )$ and shaded in the disc $\mu_{\theta}$ in Figure 2.16. 
If we
think of the subsurface $\Delta$ of Figure 2.16  as a slit-open boxing glove, the discs
$\mu_{\theta}$ for $\theta_1 < \theta < \theta_5$ will appear to be "outside" 
the glove,
whereas those for  $\theta_5 < \theta < \theta_1$ will be "inside" the glove. 
The
disc $\delta$, which is on the glove, is a limiting position for both families 
of discs.
Thus the closure of the union of all of the discs $\mu_{\theta}$, 
$\theta \in [0,2\pi]$
will be two 3-balls $B_1$ and $B_2$, which intersect along a single arc 
in the disc
$\delta$ of Figure  2.13(b), i.e. the arc which runs from $q$ to $s$ through $v$.
The arc $b(\theta_7)$  in Figure 2.17 is in the boundary of $B_1 \cup B_2$.

$$\EPSFbox{2.16-2.17.eps }$$
\centerline{\bf Figure 2.16 \ \ \ \ \ \ \ \ \ \ \ \ \ \ \ \ \ \ \ \ \ Figure 2.17}

We may assume that $B_1 \cup B_2$ intersects the link $L$ in two 
unbraided weighted arcs.
This is clearly possible, since all braiding may be pushed out of 
$B_1 \cup B_2$. Notice
that the surface $\cS$ may meet $B_1$ and $B_2$ in some number of 
sheets which are locally
parallel to the embedded disc $\Delta$ which we depicted in Figure 2.16. 
In particular,
$\cS$ may meet a neighborhood $N_A$ of $A$ in some number of radially 
foliated discs between
the discs at $w_1$ and $w_2$.

Now we need to thicken $B_1 \cup B_2$ a little bit. With this in mind, 
reparametrize the
interval $[0,2\pi ]$ so that the singularities occur at $\theta_1 = 0$ and 
$\theta_5 =
\pi$. We can then pair the discs $\mu_{\theta}$ and $\mu_{\theta + \pi}$, so that 
$\mu_{\theta}$ is in $B_1$ and $\mu_{\theta + \pi}$ is in $B_2$. 
See Figure 2.18, which
depicts subsets of $H_{\theta} \cup H_{\theta + \pi}$. Now, 
letting $0 < \theta < \pi$, 
notice (see Figure  2.13(b)) that if $H_{\theta}$ is non-singular, then $H_{\theta}$
contains two leaves in the foliation of $\Delta$: the leaf $b( \theta )$, 
which joins
$v$ to $w_2$, and also a small arc $\alpha ( \theta )$ which runs out from 
$w_1$ to the
boundary of $\Delta$. Similarly, any non-singular fiber $H_{\theta + \pi}$ 
contains the
leaf $b(\theta + \pi)$, which  joins $v$ to $w_1$, and a small arc 
$\alpha ( \theta +
\pi)$  which runs out to the boundary of $\Delta$ from $w_2$.  Let 
$N( \theta )$ be a
neighborhood of $\mu_{\theta} \cup \mu_{\theta + \pi} \cup  \alpha( \theta ) \cup 
\alpha(\theta + \pi)$ in $H_{\theta} \cup  H_{\theta + \pi}$, 
chosen so that $N( \theta)
\cap L = (\mu_{\theta}  \cup \mu_{\theta + \pi}) \cap L$ and 
$N( \theta ) \cap \cS$ is the
union of $(\mu_{\theta} \cup \mu_{\theta + \pi} ) \cap \cS$ and 
$(H_{\theta} \cup  H_{\theta
+ \pi}) \cap N_A$. Finally, choose the neighborhoods $N( \theta )$ 
so that they vary
smoothly as $\theta$ is varied between $0$ and $\pi$. Let $N$ be 
the closure of the union
of all of the $N( \theta )$, $\theta \in [0, \pi )$. Notice that 
$N$ is not a 3-ball (it
is a 3-ball with holes), because $N( \theta )$ is not a disc for 
$\theta$-values close to
$\theta =0$. Choose $\epsilon > 0$ so that $N( \theta )$ is a disc 
for $\theta \in
[\epsilon ,\pi - \epsilon ]$.

$$\EPSFbox{2.18.eps }$$
\centerline{\bf Figure 2.18}

\

We can now describe precisely the exchange move which was indicated 
earlier in 
Figure 2.15. See the top left picture in Figure 2.19, which is 
intended to correspond
to the left picture in Figure 2.18, at some $\theta \in (\epsilon ,\pi
- \epsilon )$. The plane of the paper is divided by $A$ into two
half-planes, and we will assume the left one to be $H_{\theta}$ 
and the right one to be
$H_{\theta + \pi}$. We want to describe an isotopic deformation of 
$\cS \cup L$, and shall
do so by the series of pictures in Figure 2.19, all of which 
correspond to the same
fixed value of $\theta \in (\epsilon ,\pi
- \epsilon )$.

Choose a subdisc $d$ of $\mu_{\theta + \pi}$ such that $d \cap L = \mu_{\theta +
\pi} \cap L$  and also   $d \cap \cS = \mu_{\theta + \pi} \cap \cS  \setminus    b( \theta + \pi )$.
Our isotopy is to be supported in $N(\theta )$ and is to be the identity on
$\mu_{\theta}$ and on $\mu_{\theta + \pi}$ minus a neighborhood of $d$ on $\mu_{\theta
+ \pi}$. The isotopy pushes the disc $d$ (and the points of $L$ and $\cS$ which meet it)
across $A$ and then down and eventually into the cross-hatched area below $w_2$ which was
illustrated in the blow-up in Figure 2.18. The isotopy is to be defined on pairs of
fibers, and is to be defined so that it varies continuously as we vary $\theta$. At the
end of the isotopy the link $L$ is to encircle the axis below $w_2$ instead of between
$w_1$ and $v$. Also, each sheet of $\cS$ which intersected $A$ between $w_1$ and $v$ is to
intersect $A$ the same number of times below $w_2$. (One can think of
this change as accomplished by "putting your hand into your pocket and emptying it".)

The only remaining problem is to ask what happens in the 
$\theta$-interval $[0,
\epsilon ]$? Let $u$ be the isotopy parameter. The first thing to 
notice is that when $u$
is in the interval $[.25,.75]$ the isotopy is supported in the left 
half-plane. We may
then assume that the deformed discs $d_u$, $u \in [.25,.75]$, have 
empty intersection with
$H_{\theta + \pi}$ if $\theta$ is not in the interval $[\epsilon ,\pi - \epsilon]$. 
The
second thing to notice is that when $u \in [0,.25]$ and $u \in [.75,1]$ the 
isotopy is
supported in a neighborhood of the axis $A$. Since that neighborhood can be 
chosen to be
disjoint from a neighborhood of the singular points,  it follows that
Figure 2.19 actually tells us everything we need to describe the isotopy of 
$\cS \cup L$ completely. 

At the end of the isotopy, we may reposition $\cS$ so that its 
decomposition is
exactly as it was before the change. There has, however, been a change in the
order of
the points of $A \cap \cS$ on $A$. After the change $b(\theta_7)$ 
and possibly some
leaves parallel to it will become inessential.  We can therefore modify the 
surface by an
additional isotopy to a new surface $\cS^{'}$, by eliminating all inessential 
$b$-arcs and
any resulting  unwanted  $c$-circles, in the manner described in the proofs of 
Lemmas
1.2, 1.3 and 1.4. Thus, the  decomposition of  $\cS^{'}$  has at least  two fewer
vertices and at least two fewer  regions than the decomposition of $\cS$. $\|$

$$\EPSFbox{2.19.eps }$$
\centerline{\bf Figure 2.19}

\

\

\noindent {\bf $\S$ 3. \ GLOBAL TOPOLOGY OF THE TILING OF $\cF$}. 

\

In the previous section we discussed particular ways in which we could
make local changes in the  tiling of  $\cF$ and in the decomposition of 
$\cC$. In this
section we discuss the $\lq$global
topology' of the tiling of $\cF$.  Clues about the
global topology of $\cF$ may be found in four graphs defined by the 
tiling, which we will
introduce shortly.

The following observation will be useful in several places in this section:

\

\noindent {\bf Observation $\star$}: If we change the orientation on $L$ 
then every
$\pm$ vertex becomes a $\mp$ vertex. If we change the orientation on the 
fibers of 
$H$ then every $\pm$ singular point becomes a $\mp$ point.  In this way 
certain proofs which 
assume a particular choice of sign $\epsilon$ for a vertex and sign 
$\delta$ for a singular 
point may be adapted to apply to the four possible choices of 
$(\epsilon,\delta )$. 

Our initial goal
is to locate the places on the tiling of $\cF$ where the changes of 
Section 2
can
be made.  We begin with some new definitions. Let $v$ be a vertex. A vertex $v_i$ 
in the tiling
of
$\cF$ is said to be {\it near v} if there is a $b$-arc $\beta_i$ in
the
foliation of $\cF$ with $\partial \beta_i = v \cup v_i$. A singular
point
$s_j$ in the tiling of $\cF$ is said to be {\it near  v} if $s_j$
lies
on a
singular leaf $l_j$ which ends at $v$. Let $star (v) = \{ v_i, s_j | v_i,
s_j \
are \ near \ v \}.$

See Figure 3.1 for two examples. In the first
$star(v)$ has arbitrarily many vertices, but in the second there are only two
vertices in $star(v)$.  A vertex $v$ in the tiling of $\cF$ is said to be 
an 
{\it interior} vertex
if it is not an endpoint of an $a$-arc.

$$\EPSFbox{3.1.eps }$$
\centerline{\bf Figure 3.1}

\

Our first result in this section  was  explained
to us by Menasco. He says that it is implicit in \cite{Bennequin}. 
However, it is very
difficult  to pinpoint it
there and at this writing Menasco's version of it is unpublished, so it 
appears here for the
first time.

\

\noindent {\bf Lemma 3.1} {\it Let $v$ be an interior vertex. Then $star (v)$
contains both positive and negative singular points.}

\

\noindent {\it Proof of Lemma 3.1}. By Observation $\star$  we may assume 
without loss of
generality that $v$ is a positive vertex.  Let $v_1,v_2,\dots,v_r$ be the 
vertices in $star (v)$,
ordered so that their cyclic order on the oriented braid  axis $A$ is    
$v,v_1,v_2,\dots,v_r$
(see Figure 3.2(a)). Thus there are $r$ vertices in $star(v)$, also 
\newpage
$v_1$ 
is the first vertex in
$star (v)$ which is encountered if one starts at $v$ and travels along 
the oriented braid axis
$A$ in the positive direction. 

The vertices
$v_1,\dots,v_r$ also have a second order, in the flow around $v$. 
Traveling around $v$
counterclockwise, in the direction of increasing
polar angle, let $v_0 \in \{ v_1,\dots,v_r \}$ be the vertex which occurs 
just before
$v_1$, and let $s$ be the singular point in $star (v)$ which is between 
$v_0$ and
$v_1$ (see Figure 3.2(b)).  

Choose a non-singular fiber $H_{\theta_1}$ which contains a $b$-arc
$\beta_1$ joining $v$
to $v_1$, as in  Figure 3.2 (c). The arc $\beta_1$ divides
$H_{\theta_1}$
into into two half-discs, which we call $H_{\theta_1}^+$ and
$H_{\theta_1}^-$, where the signs are
chosen so that $H_{\theta_1}^+$ (resp. $H_{\theta_1}^-$) is the side of
$H_{\theta_1}$ split along $\beta_1$
which meets the $+$ (resp. $-$) side of  $\cF$. The fact that the order of the
vertices on $A$ is $v,v_1,\dots,v_r$ then implies that the vertices
$v_2,\dots,v_r$ are all on the $\partial H_{\theta_1}^-$ part of $A$.

$$\EPSFbox{3.2.eps }$$
\centerline{\bf Figure 3.2}

\ 

Notice that the sign of a singularity has natural meaning when it is
examined in a sequence of fibers of $H$. The surface $\cF$ is  oriented, so
that
when we view a component of $\cF \cap H_\theta$ on $H_\theta$ it will 
have a
well-defined positive and negative side. See Figure 1.11 for a series of
pictures
which show the leaves of $\cF \cap H_\theta$ just before, at the instant
of, and just after a positive singularity. When
the singularity is positive (resp. negative) the negative (resp. positive)
sides
of $\cF \cap H_\theta$ approach one-another just before the singularity and
then split apart in a new way.  This means that, following the flow around
$v$, the
first singularity  occurring after $\beta_1$ must be positive. Otherwise, 
$\beta_1$ would surger  with a $b$-arc having an endpoint on 
$\partial H_{\theta_1}^+$, but  since $v_2,\dots,v_r$ are on 
$\partial H_{\theta_1}^-$, this is impossible.

Recall how the vertex $v_0$ and the singular point $s$ were chosen. In
particular, recall that $s$ is the singular point in $star(v)$ which
occurs
between $v_0$ and $v_1$. Thus, if $\beta_1$ and $\beta_0$ are the b-arcs 
joining $v$ to $v_1$
and $v_0$, respectively, then  $\beta_0,s, \beta_1$ have that order in
the fibration. When viewed on a sequence of fibers of $H$, we see that the
singularity at $s$ results in the creation of the $b$-arc $\beta_1$, and 
that just before this
singularity occurs there is a non-singular fiber $H_{\theta_0}$  
containing the $b$-arc
$\beta_0$ (see Figures 3.2 (c) and (d)). The fact that $v$ is positive 
forces $v_0$ to be
negative, because $v$ and $v_0$ are the two endpoints of a $b$-arc. This
means
that the oriented axis $A$ intersects the positive side of $\cF$ first at
$v_0$.
It follows that $v_1$ must be on the negative side of $H_{\theta_0}$ split
along
$\beta_0$. But then, the singularity at $s$ must be negative. $\|$
 
\

Let  ${\vec x} =
x, \dots,x$ where $x$ is either a sign $\pm$ or a $b$-arc.

\

\noindent {\bf Lemma 3.2} {\it After some number of changes in
foliation,  exchange moves and isotopies in the complement of the axis, 
it may be assumed
that no interior vertex   has sign $({\vec +}, {\vec -})$.} 

\

\noindent {\it Proof of Lemma  3.2}.  Let $V(\cF)$ be the set of all 
vertices in the tiling
of $\cF$.  Recall that $V$ is the cardinality of  $V(\cF)$. Suppose that
$v\in V(\cF)$ is an interior vertex  with  sign  $({\vec +}, {\vec -})$.  Since $v$ is
an interior vertex, it has type $({\vec b})$.   By Theorem 2.1  we may
assume, after changes in foliation which push off adjacent tiles of like 
sign, that
$v$ has valence 2, type  $(b,b)$ and  sign $(+,-)$.  The vertex $v$ is 
still an interior
vertex, so Theorem 2.2, part (2.2.2) applies.  The link $L$  admits an
exchange move, and after the exchange move we may eliminate an 
inessential $b$ arc which has $v$
as its endpoint, by an isotopy of $\cF$ which is supported in a 
neighborhood of the tiles which
meet $v$. Since this process eliminates $v$ and creates no new vertices, 
$V$ is reduced. If there is again an interior vertex $v$ with sign 
$({\vec +}, {\vec -})$  the argument may be repeated and,  since $V$    is reduced
each time,  a finite number of  repetitions  will yield
a tiling in which no interior vertex has sign  $({\vec +}, {\vec -})$. $\|$

\

\noindent {\bf Lemma 3.3} {\it After  some number of changes in
foliation,  exchange moves,  isotopies in the complement of the axis,  
and destabilizations
along vertices of type (a), it may be
assumed that if a vertex  $v$ has type $({\vec b},a,{\vec b})$, 
then there are
singularities of opposite sign in $star(v)$.}

\

\noindent {\it Proof of Lemma 3.3}.  Suppose there is such a vertex $v$ 
in the set  $V(\cF)$
of all vertices in the tiling of $\cF$.   By Theorem 2.1 and Remark 2.1 we may
assume, after changes in foliation which push off  tiles adjacent at a 
$b$-arc of like sign, that
$v$ has type  $(a)$. A destabilization move via the vertex $v$ then 
eliminates it. Since this
process creates no new vertices, $V$ is reduced. 
If there is
again a vertex  as described in Lemma 3.3  the argument may be repeated 
and, since
$V$ is reduced each time, we may
assume that Lemma 3.3 holds after a finite number of repetitions. $\|$

\

\noindent {\bf Lemma 3.4} {\it After some number of changes in
foliation,  exchange moves,  isotopies in the complement of the axis,  
and destabilizations
along vertices of type (a), it may be
assumed that no  vertex $v$  has   type $({\vec b},a,{\vec b})$ 
and  sign  $({\vec +}, {\vec -})$, where one of the changes from $+$ to $-$ corresponds to 
singularities in $star(v)$ of
opposite sign in $ab$-tiles adjacent at the boundary $a$-arc which meets 
$v$.}

\

\noindent {\it Proof of Lemma 3.4}.  Suppose there is such a vertex $v$ 
in the set  $V(\cF)$
of all vertices in the tiling of $\cF$.   By Theorem 2.1  we may
assume, after changes in foliation which push off  tiles adjacent at a 
$b$-arc of like sign, that
$v$ has valence 2, type  $(a,b)$ and  sign $(+,-)$. By  Theorem 2.2, part 
(2.2.1), the
link $L$  admits an exchange move and  an isotopy of $\cF$ eliminates 
$v$. Since this process
creates no new vertices, $V$ is reduced. If 
there is
again a vertex  as described in Lemma 3.4  the argument may be repeated 
and, since
$V$ is reduced each time, we may
assume that Lemma 3.4 holds after a finite number of repetitions. $\|$

\    

To continue our work, we introduce four graphs, denoted $G_{++}$, $G_{+-}$, 
$G_{-+}$ and $G_{--}$,
whose edges are singular leaves in the positive and negative 
tiles in the tiling of $\cF$. See Figure 3.3.  These graphs were first considered 
by Bennequin in \cite
{Bennequin}, who was interested in them in the special case when $\cF$ is a
disc. They are in some ways dual to the tilings, and they encode 
topological properties of
the surface in new ways. They were  used by Menasco in \cite{Menasco} in 
his study of
unknotting numbers of knots and by Fung in his thesis \cite{Fung}. No doubt
they
have other applications too. All four graphs are subsets of

$$\EPSFbox{3.3.eps }$$
\centerline{\bf Figure 3.3}

\

\noindent the singular leaves in
the foliation of $\cF$. The graph $G_{\epsilon , \delta}$ (consult Figure 3.3) 
passes through  only 
vertices of sign $\epsilon$ and singular points of sign $\delta$,
where $\epsilon$ and $\delta$
are fixed signs, each $\pm$. However, they have rather different 
structures and
pick up different aspects of the combinatorics of the tiling, because the
three
tile types in a spanning surface are quite different as regards their
positive
and negative vertices.  Recall that the {\em sign} of a tile is the sign 
of the unique singular
point in that tile.

See
Figure 3.3 for illustrations of the edges and vertices in $G_{++}$, 
$G_{+-}$, $G_{-+}$ and
$G_{--}$. The black vertices indicate vertices of $G_{+,\delta}$ and the 
circled vertices indicate
vertices of $G_{-,\delta}$. The thick black lines indicate edges of 
$G_{\epsilon,+}$, and the double
lines indicate edges of $G_{\epsilon,-}$. Notice that the graphs 
$G{+,-}$ 
and $G_{-,+}$ are obtained from $G_{+,+}$ and $G_{-,-}$ by reversing the orientation of 
the fibers of $H$. 

\

\noindent {\bf Definition} (See Figure 3.3): 
\begin{itemize} 
\item  The graph $G_{+,\delta}$  has as its edges those subarcs of 
singular leaves which
join the two + vertices in a $aa$,$ab$ or $bb$ tile of sign $\delta$. The 
vertices
of  $G_{+,\delta}$ are the endpoints of these edges, together with all + 
vertices in the tiling
of $\cF$ which are not adjacent to any singular point of sign $\delta$.  

\item  The graph $G_{-,\delta}$ has as its edges those subarcs of 
singular leaves which:
\begin{itemize}

\item join the two negative vertices in a $bb$ tile of sign $\delta$.

\item join the negative vertex in an $ab$ tile of sign $\delta$ to 
$\partial\cF$.

\item join $\partial\cF$ to $\partial\cF$ in an $aa$ tile of sign $\delta$.
\end{itemize}  

The vertices of $G_{-,\delta}$ are the endpoints of 
these edges, together with negative
vertices in the tiling of $\cF$ which are not adjacent to any singular 
point of
sign $\delta$.
\end{itemize}

\noindent{\bf Remark 3.1} If one thinks of $\partial \cF$ as acting 
"like a negative vertex" (see Remark 1.1) then the
definition of $G_{-,\delta}$ will be seen to be analogous to that of 
$G_{+,\delta}$. 

\

\noindent{\bf Remark 3.2} Note that there are vertices of 
$G_{\epsilon,\delta}$ that are not
vertices of tiles in the tiling of $\cF$. To avoid confusion we will 
refer to
the former as {\it vertices of $G_{\epsilon,\delta}$} and the later as 
{\it vertices in the
tiling of $\cF$}.

\

\noindent {\bf Lemma 3.5} \cite{Menasco}. {\it The graphs 
$G_{\epsilon,\delta}$ have the
following properties}:

\begin{enumerate}

\item {\it $G_{\epsilon,\delta} \cap G_{-\epsilon, -\delta} = \emptyset$.}

\item {\it Every singular point in the foliation of $\cF$ is
in $G_{++}$ or $G_{--}$ (and so also in $G_{+-}$ or $G_{-+}$).} 

\item {\it Every  vertex in the tiling of $\cF$ is a vertex of $G_{++}$
or $G_{--}$ ( and so also in $G_{+-}$ or $G_{-+}$).}

\end{enumerate}

\noindent {\it Proof of Lemma 3.5}: It suffices to prove the assertions for
$G_{+,+}$ and $G_{-,-}$, since the other two cases follow by reversing 
the orientation of fibers
of $H$. 
\begin{enumerate}
\item The vertices and singular points of $G_{++}$ and $G_{--}$ have
distinct signs. Since each edge of $G_{++}$ and $G_{--}$ contains a 
singular point, the
only possibility is that $G_{++} \cap G_{--} = \emptyset.$ 

\item Each singular point $s$ in the foliation of $\cF$ belongs to either an
$aa$ or $ab$ or $bb$-tile, and a check of the possibilities shows that $s$
belongs to an edge of $G_{++}$ if the sign of $s$ is
positive and an edge of $G_{--}$ if its sign is negative.

\item This follows from (2). $\|$
\end{enumerate}

\ 

A vertex of
$G_{\epsilon,\delta}$ is {\it isolated} if it is not the endpoint of any 
edge of
$G_{\epsilon,\delta}$. 

\

\noindent {\bf Lemma 3.6}. {\it No graph $G_{\epsilon,\delta}$  contains 
an isolated
vertex which is an  interior vertex in the
tiling of $\cF$.} 

\

\noindent {\it Proof of Lemma 3.6} Suppose that for some choice of 
$(\epsilon,\delta)$ the
vertex $v$ is an isolated interior vertex in $G_{\epsilon,\delta}$. Since 
$v$ is an isolated
vertex, every singular point in $star(v)$ must have sign $-\delta$.  
However this contradicts 
Lemma 3.1.  $||$

\

An {\it endpoint} of  $G_{\epsilon,\delta}$ is  a vertex of 
$G_{\epsilon,\delta}$ which is the endpoint of exactly one edge of 
$G_{\epsilon,\delta}$.

\

\noindent {\bf Lemma 3.7}  {\it After some number of changes in
foliation,  exchange moves and isotopies in the complement of the axis, 
it may be assumed
that no graph  $G_{\epsilon,\delta}$ has an endpoint vertex which is an  
interior vertex in
the tiling of $\cF$.} 

\

\noindent {\it Proof of Lemma 3.7}. Suppose   $v$ is an endpoint vertex of
$G_{\epsilon,\delta}$ and an interior vertex in the tiling of $\cF$, for 
some fixed $\epsilon$
and $\delta$. Since $v$ is an interior vertex in the tiling of $\cF$ we 
may assume by Lemma  3.2, 
after some number of changes in foliation,  exchange moves and isotopies 
in the complement 
of the axis, that  $v$  does not have  sign $({\vec +}, {\vec -})$. 
But,
since $v$ is an endpoint vertex of $G_{\epsilon,\delta}$,  there is
exactly one singularity in $star(v)$ with sign $\delta$ and all others 
have sign $-\delta$. Thus,
$v$ has sign $(\delta, {\vec -\delta})$, a contradiction. $\|$

\

\noindent {\bf Lemma 3.8} {\it After some number of changes in foliation, 
exchange moves, 
isotopies in the complement of the axis, we may assume that:

\begin{enumerate}

\item[(i)] no closed loop  in $G_{\epsilon , \delta}$   bounds a disk on 
$\cF$;

\item[(ii)] no closed loop which is the union of an edgepath in 
$G_{\epsilon, +}$ and    an
edgepath in $G_{\epsilon, -}$ bounds a disk on $\cF$;

\item[(iii)] after some number of  destabilizations along vertices of 
type (a), no closed loop which  is a
union of edges in $G_{-, \delta}$ and one subarc of $L$   bounds a disk 
on $\cF$ and

\item[(iv)] after some number of  destabilizations along vertices of type 
(a), no closed loop which is the
union of an edgepath in $G_{-, +}$,     an edgepath in
$G_{-, -}$ and a subarc of $L$ bounds a disk on $\cF$;

\end{enumerate}}

\noindent {\it Proof of Lemma 3.8}. By Lemmas  3.2, 3.6 and 3.7, we may 
assume, after  some
number of changes in foliation, exchange moves,  isotopies in the 
complement of the axis, that
no graph $G_{\epsilon, \delta}$ contains  an interior vertex in the 
tiling of $\cF$
which has  sign  $({\vec +}, {\vec -})$, and no graph 
$G_{\epsilon, \delta}$ contains  an
isolated  or endpoint vertex that is also an interior vertex in the 
tiling of $\cF$. 

Suppose (i) or (ii) is false. Then there is a closed loop   $c$  of the type
described in (i) or (ii) which bounds a disk $D$ on $\cF$. Choose $c$ 
to be  minimal in the
sense  that there is no closed loop as  described in (i) or (ii)  which 
bounds a proper
subdisk  of $D$. The loop    $c$ is either  in
$G_{\epsilon,\delta}$ for some choice of
$\epsilon, \delta$, or is equal to the union of an edgepath in 
$G_{\epsilon,+}$
and an edgepath in $G_{\epsilon,-}$, for some choice of $\epsilon$.  
Therefore  any edge $e$ of $c$ is in  $G_{\epsilon,\mu}$, where $\mu 
= 
1$ or $\mu = -1$. Thus,  since
$D \cap L = \emptyset$,  $e$ is in a tile which  contains a  vertex $w$ 
in $D$
of  parity $-\epsilon$ as shown in  (3),(5) or (6) of Figures 3.3(a)(b),
depending on the signs of $\epsilon$, $\mu$. Since
$w$ has parity $-\epsilon$, $w$ is in $Int(D)$. 

We have just shown that $G_{-\epsilon, \mu}$ has non-empty
intersection with $Int(D)$ for $\mu - \pm 1$, but what about $G_{\epsilon, \mu}$?
Note that   $G_{\epsilon,\mu} \cap Int(D)$ contains only tree components with
endpoints on $c$, since otherwise  either $c$ is not minimal or
$G_{\epsilon,\mu}$ contains an  endpoint or isolated vertex which is 
an interior vertex in the tiling of $\cF$. Suppose $T$ is such a tree 
component. Note that each endpoint of $T$ must be a vertex of $c$.
For, if some endpoint of $T$ is a singularity of $c$,  that singularity 
must be in a tile having
three vertices of the same parity $\epsilon$. But since no tile type has 
more than  two vertices
of the same parity, we conclude that $T$ has an edgepath $\alpha$ with 
endpoints equal to
vertices of $c$. The edgepath  $\alpha$ separates $D$ into two subdisks. 
But each subdisk 
has boundary in $G_{\epsilon,\mu}$   or has boundary equal to the union of
an edgepath in $G_{\epsilon, +}$ and    an edgepath in $G_{\epsilon, 
-}$,  contradicting 
the minimality of $c$. We may therefore assume that $G_{\epsilon,\mu} 
\cap Int(D) = \emptyset$
for $\mu = \pm 1$. 

By Lemma 3.5 (2),  no singularities in the 
foliation of $\cF$ are
contained in $Int(D)$. But then   all singularities in $star(w)$  are on 
$c$ and each  vertex of
$c$ is connected to  $w$  by a $b$-arc. Thus, $w$  has  type $({\vec b})$,
and sign $({\vec \delta})$, if $c$ is in $G_{\epsilon, \delta}$ or
sign $({\vec +}, {\vec -})$ if $c$ is the union of an edgepath in
$G_{\epsilon, +}$ and $G_{\epsilon, -}$. The first case contradicts Lemma 3.1
and the second contradicts Lemma  3.2. This proves (i) and (ii).

We now assume that, after some number of  destabilizations along vertices of
type (a), that Lemmas 3.3 and 3.4 hold in addition to the Lemmas  3.2, 
3.6 and 3.7. Suppose
that  (iii) or (iv) is false. Then there is a closed loop   $c$  of the type
described in (iii) or (iv) which bounds a disk $D$ on $\cF$. Choose $c$ 
to be  minimal in the
sense  that there is no closed loop as  described in (iii) or (iv)  which 
bounds a proper
subdisk  of $D$.   Note that every edge $e$ of  $c$ is in  
$G_{-,\delta}$  for some choice of
$\delta = \pm 1$, so $e$ is as shown in one of Figures (2)(4)(6) of  3.3 
(a) (respectively 3.3
(b)) if $\delta = -$ (respectively  $\delta = +$). Thus, there is a 
positive vertex $w$ in $D$.
Since every vertex of $c$ is negative, $w$ is in $Int(D)$. Note that 
$G_{-,\delta} \cap Int(D) = \emptyset$  for $\delta = \pm 1$. For, 
suppose otherwise. By (i),
there are no closed loops, and  since every vertex in $G_{-,\delta} \cap 
Int(D)$ is an interior
vertex in the tiling of $\cF$, $G_{-,\delta} \cap Int(D)$ contains no 
endpoint or isolated
vertices. Thus, $G_{-,\delta} \cap Int(D)$ consists of  tree components 
with endpoints on $c$,
and $T$ contains  an edgepath $\alpha$  with endpoints on $c$. No  
endpoint of $\alpha$ can be a 
singularity of $c$, since otherwise there is a tile with three negative 
vertices. Thus, the
endpoints of $\alpha$ may be points on $L$ or vertices  of $c$. It is 
easy to
check that all possible cases contradict (i), (ii) or minimality of $c$, 
and we conclude that 
$G_{-,\delta} \cap Int(D) = \emptyset$ for  $\delta = \pm 1$. By Lemma 
3.5, we may 
assume that $Int(D)$ contains  no singularities in the foliation of 
$\cF$. Thus, all
singularities of $star(w)$ are on $c$ and each  vertex of $c$ is 
connected to  $w$  by a $b$-arc.
Thus, $w$  has  type $({\vec b}, a, {\vec b})$.  If $c$ is a loop 
of the type described
in (iii), then all singularities in $star(w)$ have the same sign, 
contradicting Lemma 3.3. If
$c$ is a loop of the type described in (iv), then $w$ has sign  $({\vec +}, {\vec
-})$, where one of the  changes from $+$ to $-$ corresponds to singularities in 
$star(w)$ of
opposite sign in $ab$-tiles adjacent at the boundary $a$-arc which meets 
$w$. But this
contradicts Lemma 3.4 and concludes the proof of Lemma 3.8. $\|$

\

\noindent {\bf Theorem 3.1.} {\it After some number of
changes in foliation,  exchange moves and isotopies in the complement of 
the axis, it may be
assumed that each of the following holds for all four graphs at once:

\begin{enumerate}

\item[(i)] $G_{\epsilon,\delta}$ has no interior endpoint vertices.

\item[(ii)] $G_{\epsilon,\delta}$ has no interior isolated vertices.

\item[(iii)] $G_{\epsilon,\delta}$ contains no closed loop  bounding a 
disk on $\cF$.

\item[(iv)] there is no closed loop equal to the union of an edgepath in 
$G_{\epsilon, +}$ and
an edgepath in $G_{\epsilon, -}$ which bounds a disk on $\cF$.

\item[(v)] after some number of destabilizations at vertices of type (a), 
there is no closed
loop equal to the union of an edge-path  in $G_{-,\delta}$ and a subarc 
of $L$ which 
bounds a disk on
$\cF$.

\item[(vi)] after some number of deletions of trivial loops of $L$, there 
is no closed loop equal
to the union of an edge-path  in $G_{-,+}$, an edge-path in $G_{-,-}$ and 
a subarc of $L$ which 
bounds a disk on
$\cF$.

\end{enumerate}}

\noindent {\it Proof of Theorem 3.1}. This follows from Lemmas 3.6, 3.7 
and 3.8. $\|$

\

\

\noindent {\bf $\S$ 4. APPLICATIONS AND OPEN PROBLEMS}. 

\

In this section we
prove three results (Theorems 4.1, 4.2 and 4.3) which
illustrate how the methods which were described in Sections 1-3 may be used to
solve problems about surfaces which have a special position with regard to a
link. We also pose some interesting open problems which should yield through the use
of  the techniques in this paper. With regard to the theorems: Theorem 4.1 is from
\cite{Birman-Menasco VI}. To the best of our knowledge Theorem 4.2 has not been
written down anywhere, although it is implicit in the work of Bennequin and
Birman-Menasco. Theorem 4.3 is the main result from \cite{Birman-Menasco V}, but
the proof which we give here is new and is simpler than the original proof,
because we use the graphs $G_{\epsilon,\delta}$, which were not known when
\cite{Birman-Menasco V} was written. 

Our first application relates to Theorem 1.2, in the special
case when $\cS = \cF$. Let us
suppose that $L$ is a closed $n$-braid with respect to the braid axis $A$ and that $\cF$ is
decomposed into tiles of types $aa$, $ab$ and $bb$. Suppose that there are $t$
singularities in the foliation of $\cF$, that $A$ intersects $\cF$ in  $k$ points, and that
there are $r$ points on $L$ which are the endpoints of singular leaves in the tiling. Assume
that the singularities $s_1, s_2,  \dots , s_t$ have been numbered in the cyclic order in
which they occur in the fibration,  the vertices have been numbered
$v_1,v_2,\dots,v_k$ in the natural order in which they occur on $A$, and  the endpoints
of the singular leaves which are on $L$ have been numbered $p_1,p_2,\dots,p_r$ in their
natural cyclic order on $L$. 

\

\noindent {\bf Theorem 4.1}. {\it The following
combinatorial data for the tiling of $\cF$ determines the embedding of $L$
(respectively $\cF$)}: 

\begin{enumerate}

\item[(1)] {\it A listing of the number of regions of type $aa$ and $ab$
(respectively $aa, ab$ and $bb$).}   

\item[(2)] {\it For each region, its
sign and an identification of its vertices among the vertices
$v_i$ on $A$.}

\item[(3)] {\it For each region of type $aa$ and $ab$, an identification of its
$L$-endpoints among the points $p_j$ on $L$.}

\item[(4)] {\it For each region,  the order of  the  singularity it
contains  in the cycle  $s_1, \dots , s_t$.}

\end{enumerate}

\noindent  {\it Proof of Theorem 4.1.} In the proof of Lemma 1.6 we showed how
each $ab$-tile has a canonical embedding in 3-space, up to an isotopy which
is supported on the axis $A$ and the singular fiber $H_\theta$ associated
to the tile. The embedding is determined by the data of type (2) above. A small
modification of the argument shows that each $aa$ and $bb$-tile also has a similar canonical
embedding. See Figure 1.13 of $\S$1 for illustrations of the canonical embeddings
of the three tile types.  The embedding of $\cF$ is determined  by specifying how the
boundary components of the embedded regions are identified. Identifications of boundary
components which are $a$-arcs are specified by the data of type (3). Two $b$-arc boundary
components of regions $R_1$ and $R_2$  are identified if they have two common
vertices and the singularities belonging to the regions $R_1$ and $R_2$ occur
consecutively in the fibration. Thus, the data of type (2) and (4) determines how the
boundary components of regions which are $b$-arcs are identified. The embedding of $\cF$ is
therefore determined by the data of types (1), (2), (3) and (4).  If we are
only interested in the embedding of $L$, then we do not need to include 
the  information for $bb$-tiles  because they do not meet $L$. $\|$  

\

\noindent  {\bf Problem 4.1}. A word in the standard elementary braid
generators $\sigma_1,\dots,\sigma_{n-1}$ of the braid group $B_n$ and their
inverses gives a description of a closed braid by a finite set of combinatorial
data. Theorem 4.1 gives an alternative description, which in principle contains
more information. In principle, the new data implies the old data, however
we do not at this time know any direct algorithm which allows us to pass from the
combinatorial data of Theorem 4.1 to a braid word which describes $L$, although
it is clear that such an algorithm must exist. It would be
interesting to have such an algorithm. The key to the difficulty seems to be to
understand the way in which $ab$-tiles affect the braid word. 

\
         
\noindent  {\bf Remark 4.1}: Theorem 4.1 generalizes to the
case of closed surfaces, if one includes additional data on the ways in which
regions of type $bb$, $bc$ and $cc$ are pasted together along their boundaries.

\

\noindent  {\bf Problem 4.2}. We showed, in Theorem 4.1, that the tiling
(plus some additional combinatorial data)  determines the embedded
surface.  On the other hand, it is fairly simple to construct a tiling of a
surface of genus $g$ with boundary which is not geometrically realizable by a
spanning surface for a closed braid. We pose, as an open problem, to determine
necessary and sufficient conditions such that an abstract decomposition of a
surface into foliated regions, as described in Section 1, is realized by an
embedded surface. There are obvious variations on this: e.g. when is the surface
incompressible etc. 

\ 

Recall that a
Markov surface is, by definition, a surface of minimal genus with boundary a
given closed braid. The closed braid will not, in general, have minimal braid
index. A {\it special Markov surface} is a Markov surface which is tiled entirely
by $aa$-tiles. Our next result shows that every link is the boundary of a
special Markov surface. 

\

\noindent  {\bf Theorem 4.2} {\it Every link type $\cal {L}$ has a closed
braid representative $L^*$ which is is the boundary of a special Markov surface
$\cF^*$.}

\

\noindent  {\it Proof of Theorem 4.2}. Choose a Markov surface $\cF$ for any closed braid
representative $L$ of $\cal {L}$. If the tiling of $\cF$ contains $ab$-tiles, we may stabilize
along them, as in $\S$2.2.1. Each stabilization move deletes a negative vertex at the expense
of adding a trivial loop. Thus, although new $ab$-tiles may be created from $bb$-tiles during
a stabilization move, finitely many stabilizations will eliminate all of the 
$ab$-tiles, since $\cF$ has finitely many negative vertices. Let $\cF^*$ be the
modified surface with no $ab$-tiles. Suppose the tiling of  $\cF^*$ contains a
$bb$-tile.  This $bb$-tile must be contained in a  component
of $\cF^*$ tiled entirely by $bb$-tiles, since any connected component containing both
$aa$-tiles and $bb$-tiles must contain at least one $ab$-tile. But then $\cF^*$ contains a
closed component, a contradiction. Therefore, $\cF^*$ is tiled entirely by $aa$-tiles, i.e.
$\cF^*$ is a special Markov surface. $\|$    

\

\noindent {\bf Corollary to Theorem 4.2  (Markov's Theorem for the case of the
unknot) } \cite{Bennequin}. {\it Let $K$ be any closed braid representative of the unknot.
Let $U$ be the standard 1-braid representative. Then there exists a finite sequence of
closed braid representatives $K = K_0 \to K_1 \to \dots K_r = U$ and 
$1 \leq p \leq r$ such that if $i \leq p$ (resp. $i>p$)  then $K_i$ is obtained from
$K_{i-1}$ by a single stabilization along an $ab$ tile (resp. destabilization via a type
$(a)$ vertex.}    

\

\noindent {\it Proof of the Corollary:}  \ Since $K$ is the unknot, it is the boundary of a
Markov surface $\cF=\cF_0$, which is a topological disc.  By Theorem 2.2 we may stabilize
along $ab$ tiles to obtain from $\cF_0$ a special Markov surface which we call  $\cF_p$.
Its boundary is $K_p$, and the moves from $K_0$ to $K_p$ are as claimed. This new surface is
also a disc, and it is tiled entirely by $aa$ tiles. In the new tiling the graphs $G_{+,+}$
and $G_{+,-}$ must both be trees, in view of Theorem 3.1 and the fact that $\cF_p$ is a disc.
Each endpoint vertex is then a vertex of type $(a)$. Destabilizing, we reduce the number of
vertices but leave the tiling unaltered outside the tile we just removed. In this way the
endpoint vertices may be removed, one after another, until we obtain a special Markov
surface $\cF_r$ which is  radially foliated. Its boundary $K_r$ is the standard 1-braid
representative $U$ of the unknot. $\|$ 

\

\noindent {\bf Problem 4.3}. A special Markov surface is said to be a {\it
Bennequin surface} if its boundary has minimal braid index. In the manuscript
\cite{Birman-Menasco III} it is proved that every link which is represented as a
closed n-braid with $n \leq 3$ is the boundary of a Bennequin surface. We ask: what are
necessary and sufficient conditions for a link to bound a Bennequin surface.  

\

Our final application is to give a simple proof of the main result from
\cite{Birman-Menasco V}, i.e. that there is a systematic procedure for
simplifying an arbitrary closed braid representative of the $r$-component unlink to
the trivial closed $r$-braid.   The
simple proof uses the results in $\S 3$, which were
not available when this result was first proved. After we complete the proof we will
illustrate its meaning via an example:

\

\noindent  {\bf Theorem 4.3  (Markov's theorem without stabilization in the case of the unlink)}
\cite{Birman-Menasco V}. {\it Every closed braid representative $L$ of the $r$-component unlink
${\cal U}_r$ may be reduced to the standard $r$-braid representative $U_r$ by a finite sequence of
the following moves: exchange moves, isotopy of $L$ in the complement of the braid
axis, and destabilizations along vertices of type $(a)$.}

\

\noindent  {\it Proof of Theorem 4.3}. Choose an arbitrary closed braid
representative $L$ of ${\cal U}_r$. Our closed braid
$L$ is the boundary of a Markov surface $\cF$ which is a union of $r$
disjoint discs $D_1,\dots,D_r$, each of which admits a tiling. We will show
that,   after some
number of changes in foliation, exchange moves, isotopies in the complement
of the axis and
destabilizations along vertices of type $(a)$, each $D_i$ is foliated entirely by
$a$-arcs.

By Theorem  3.1 (iii), we may
assume that no closed loop in $G_{\epsilon,\delta}$ bounds a disk on any disk $D_i$.  Since
every closed loop on a disc bounds a disc, we conclude  that each component of
$G_{\epsilon,\delta}$ is a tree or an isolated vertex. By Theorem 3.1 (ii), each component
of $G_{-,-}$ is a tree. Suppose there is such a tree component $T$. Let $\alpha$ be an
edgepath of $T$ connecting two endpoints of $T$. By Theorem 3.1 (i), $\partial \alpha
\subset L$. But then $\alpha$ separates some $D_i$ into two subdisks, each of which is
bounded by loop consisting of edges in $G_{-,-}$ and a subarc of $L$. Since this contradicts
Theorem 3.1 (v), we may assume that $G_{-,-} = \emptyset$. 

Suppose there is a tree component $T$ of $G_{+,+}$. If  $v$ is an  endpoint
vertex of $T$,  then $star(v)$ contains exactly one positive vertex. Since $G_{-,-} =
\emptyset$, there are no negative singularities in $star(v)$. Thus, $v$ has valence 1,  and
it may be eliminated by a destabilization move. Induction on the number
of endpoints of $T$, allows us to conclude that $G_{+,+}$ consists entirely of isolated
vertices. But then there are no singularities in the foliation of $\cF$ and, consequently,
each $D_i$ is foliated entirely by $a$-arcs. $\|$

\

\noindent {\bf Example 4.1.} We illustrate Theorem 4.3 with an example. 
The example, which is indicated  in Figures 4.1, suggests  an infinite sequence
of closed $4$-braids, all of which determine the same knot type. We examine this
example in the special case of the unknot which occurs when  $X=\sigma_2$ and
$Y = \sigma_2\sigma_1$.

$$\EPSFbox{4.1.eps }$$
\centerline{\bf Figure 4.1}

\

\noindent Thomas Fiedler \cite{Fiedler} has introduced an
invariant of conjugacy class and used it to prove that when  $X=\sigma_1\sigma_2$ and
$Y = \sigma_2$  these braids are in infinitely many distinct conjugacy classes, that
is, no two closed braids in the sequence are isotopic in the complement of the braid
axis $A$. We will see in our case, which is similar,  that the foliation of the disk bounded by
the unknot indicates that   any two adjacent braids in the sequence are related by an
exchange move. (Remark: since
exchange moves {\em reduce} complexity, the arrows $4.1(c) \to 4.1(b) \to 4.1(a)$  are
in the correct direction). The foliation of the disk bounded by the unknot in Figure
4.1(a)  will then show  that the
standard representative of the unknot can be obtained after destabilizations along
vertices of type $(a)$. Thus, this example illustrates how a very complicated
braid representative for the unknot may be transformed into the standard one by a
sequence of the moves described in  Theorem 4.3.

To begin to understand this example in terms of the foliated disks which the
knot bounds, the reader should first verify that the knot projections shown in
Figures 4.1(c),(b),(a) correspond to the  knots embedded in $3$-space of Figures
4.2(c),(b),(a), respectively.   The tilings of
the disks bounded by these knots are shown in Figures 4.3(c),(b),(a), respectively. To
see   that the  embedded surfaces in Figures 4.2 are in fact discs and that the
tilings in Figures 4.3 are the correct ones for the disks in Figures 4.2, it is
helpful to first study Figures 4.2(a), 4.3(a) and then the transition from Figures
4.2(a), 4.3(a) to Figures 4.2(b), 4.3(b), respectively.

The reader should have no trouble verifying that surface in Figure 4.2(a) is  a
disk and that the singular leaves of its tiling are as shown in Figure 4.3(a). 
Note that, running along the boundary of the closed
braid in Figures 4.2(a) and 4.3(a), starting at the point $q_+$, one encounters the
points $q_+,q_-,r_-,r_+,p_+,p_-$ in that order. These points are the endpoints on the
unknot of subarcs $p_+p_-$, $q_+q_-$, and $r_+r_-$ of the three singular leaves in the
tiling shown in 4.3(a).  The points where the axis  pierces the disc in 4.2(a)  are
numbered $1,2,3,4$. The singular leaf $p_+p_-$ cuts off a subdisc $P$ which contains
$1$, and similarly there are subdiscs $Q,R$ containing $3,4$ which are cut off by
$q_+q_-$, $r_+r_-$. The center subdisc, containing $2$, will be called $W$. From
Figure 4.2(a) we see that the subdiscs $P,W,R,Q$ are pierced by $A$ in that order.
Note that the surface in Figure 4.2(b) can be obtained from the one in Figure 4.1(a) by
first creating $\lq\lq$pockets" in $P$ and $W$ (which add new vertices $u_1,u_2,v_1,v_2$ to
the tiling) and then lifting the disk $Q$ up and around and  under the braid $X$ and
then pushing it back up onto the axis again by slipping it into the pockets in $P$
and $W$. Thus, the surface in Figure 4.2(b) is a disk. To see that its tiling is as
shown in Figure 4.3(b),  note that each pocket is formed by the union of two
$ab$-tiles joined at their bounding $a$-arcs and $b$-arcs. The singular leaves of the
$ab$-tiles must be as shown in Figure 4.3(b), due to the cyclic order in the fibration
of the singularities in the tiling of  Figure 4.2(b). 

Similarly, the embedding shown in Figure 4.2(c) can be obtained from the embedding
in  Figure 4.2(b) by creating new pockets in $P$ and $W$ and slipping the old pockets
together with $Q$ around $X$ and up  into the new pockets through the axis. Again,
the tiling corresponding to Figure 4.2(c) must be as shown in Figure 4.3(c), due to
the cyclic order in the fibration
of the singularities in the tiling of  Figure 4.2(c).  Continuing in this manner,

$$\EPSFbox{4.2.eps }$$

\noindent we may obtain embeddings and tilings for all disks bounded 
by our infinite sequence of unknots.

To see that   the passages from Figures (c) to (b) to (a) illustrate the moves
of Theorem 4.3, note that each of the  vertices $v_1$, $v_2$, $v_3$, $v_4$ in Figure
4.3(b)(c) has  valence 2, type $(b,b)$ and sign $(+,-)$. Thus, by Theorem 2.2,
each such vertex determines an exchange move. The passage from Figure 4.3(c)  to (b) is
realized by the two exchange moves which  are  determined by the
vertices $v_3$, $v_4$, and  the passage from Figure 4.3(b)  to (a) is realized by
the two exchange moves which are  determined by the vertices $v_1$, $v_2$. Thus, the
moves which take us from  Figures 4.2(c) to 4.2(b) and  4.2(b) to 4.2(a) are in each
case two exchange moves, each of which  empties a pocket and is 
followed by deletion
of inessential  $b$-arcs, i.e. by the collapse of the pocket. In this manner, we can
pass from any knot in our infinite sequence down to Figure 4.2(a) via exchange
moves. Now note that the vertices $1,3,4$ in Figure 4.3(a) all have type $(a)$. Thus,
after destabilizing along each of these vertices,  we
obtain the standard representative of the unknot  from the  disk in Figure 4.2(a).

\

\noindent {\bf Remark 4.2} \ \ Example 4.1 does
not really illustrate the full power of Theorem 4.3. There are vertices of type
$(a)$ in the foliated discs of Figures 4.3(c),(b),(a) and so our foliated discs
can be simplified to yield 3-braid examples via destabilization.  Since all
closed 3 and 2-braid representatives of the unknot can be reduced to the
trivial 1-braid representative without using exchange moves \cite{Birman-Menasco
III} the need for exchange moves will thus disappear.

There is another closed 4-braid representative of the knot, due to
Hugh Morton (see
\cite{Birman-Menasco V}), which is more complicated than the one we have shown
because the associated foliated disc does not contain a vertex of type
$(a)$. By Theorem 4.3 the only way one can simplify that foliation without
increasing braid index is by first using (two) exchange moves. Unfortunately,
however, that example was simply too complicated for us to be able to illustrate
it (and the pictures in Figure 4.2 should make it clear why)!

Our examples do,
nevertheless, illustrate a related point. If one goes back to Figure 4.1
and chooses  arbitrary 3-braids $X$ and $Y$ in place of the braids
$\sigma_2$ and $\sigma_2\sigma_1$ which were illustrated in 
Figures 4.2 and 4.3,
then one may easily obtain a more complicated knot $K$ which bounds a surface
$\cF$ (now of genus $>0$), and the foliation of that surface will 
illustrate the
need for exchange moves because there will not be any vertices of type $(a)$.
On the other hand, foliated surfaces of genus $>0$ present their own
difficulties when one tries to illustrate them, and so we settled for an
 example
which illustrated the phenomenon in a case  where it was
possible to make it concrete.

$$\EPSFbox{4.3.eps }$$
\centerline{Figure 4.3}

\

\noindent {\bf Remark 4.3.} Theorem 4.3 is a much stronger result than the
Corollary to Theorem 4.2, or Markov's well-known theorem for arbitrary knots and
links, because the moves which are used in Theorem 4.3 either preserve or
reduce braid index. Thus Theorem 4.3 is $\lq\lq$Markov's
theorem without stabilization", in the special case of the unlink. 
On the other hand, Theorem 4.2 illustrates the powerful role which stabilization plays,
enabling us, as it does, to remove all the `bubbles' from an arbitrary Markov surface and
replace it with one which has minimal genus and only discs which are joined by unknotted
half-twisted bands. Unfortunately, however, these `nice' Markov surfaces are
not, in general, bounded by closed braids which have minimal braid index for the 
given knot type.

\

In the case of more general knots and links, Markov's theorem without stabilization
requires
{\it braid preserving flypes} in addition to the
exchange moves of Theorem 4.3. A rudimentary version of the braid-preserving flype plays an
important role in the work in \cite{Birman-Menasco III}, which is Markov's theorem without
stabilization in the special case of links which are closed 3-braids. 

\

\noindent {\bf Problem 4.4} In the manuscript \cite{Los} J. Los has
applied ideas from \cite{Birman-Menasco ET} to generalize Theorem 4.3 to the
case of braided satellite knots. He has proved that by repeated use of isotopy,
exchange moves and destabilization along vertices of type $(a)$  any closed braid
representative of a braided satellite knot can be reduced to minimum braid index. Such a result
is false for arbitrary closed braid representatives of arbitrary links, nevertheless there
should be very general cases under which it is true. We pose this as an open problem. 

\

\noindent {\bf Problem 4.5} In the same manuscript \cite{Los} J. Los has proved
that torus knots (like the unknot and unlink) admit a unique closed braid
representative of minimum genus. Examples are given in \cite{Birman-Menasco
III} and general mechanisms are constructed in \cite{Birman-Menasco VI} which
show that this is not true in general. We pose the open problem of finding other
cases when it is true. 

\

\noindent  {\bf Problem 4.6}  How can one generalize the graphs $G_{\epsilon,\delta}$ to
related graphs on foliated closed embedded surfaces?  Do they simplify the
results in \cite{Birman-Menasco ET}, \cite{Finkelstein},
\cite{Birman-Menasco IV}?

\end{document}